# THE ART OF FRAME THEORY

PETER G. CASAZZA

ABSTRACT. The theory of frames for a Hilbert space plays a fundamental role in signal processing, image processing, data compression, sampling theory and more, as well as being a fruitful area of research in abstract mathematics. In this "tutorial" on *abstract* frame theory, we will try to point out the major directions of research in abstract frame theory and give some sample techniques from each of the areas. We will also bring out some of the important open questions, discuss some of the limitations of the existing theory, and point to some new directions for research.

## 1. INTRODUCTION

Although the Fourier transform has been a major tool in analysis for over a century, it has a serious lacking for signal analysis in that it hides in its phases information concerning the moment of emission and duration of a signal. What was needed was a localized time-frequency representation which has this information encoded in it. In 1946 D. Gabor [84] filled this gap and formulated a fundamental approach to signal decomposition in terms of elementary signals. Gabor's approach quickly became a paradigm for the spectral analysis associated with time-frequency methods. Today, Gabor's ideas are still at the center of the myriad of applications of Gabor (Weyl-Heisenberg) frames. Gabor went on to receive the Nobel Prize in Physics in 1971 for his development of holography.

Frames for a Hilbert space were formally defined by Duffin and Schaeffer [63] in 1952 to study some deep problems in nonharmonic Fourier series. Basically, Duffin and Schaeffer abstracted the fundamental notion of Gabor for studying signal processing. The ideas of Duffin and Schaeffer did not seem to generate much general interest outisde of nonharmonic Fourier series however (see Young's book [139]) until the landmark paper of Daubechies, Grossmann and Meyer [57] in 1986. After this groundbreaking work, the theory of frames began to be more widely studied, although not to the extent of the extremely rapid development of wavelets.

Traditionally, frames have been used in signal processing, image processing, data compression and sampling theory. Today, ever more uses are being found

1991 *Mathematics Subject Classification*. Primary: 42C15, 42A38.
Ths author was supported by NSF DMS 9706108.

for the theory such as optics, filterbanks, signal detection, as well as the study of Besov spaces, in Banach space theory etc. In the other direction, powerful tools from operator theory and Banach space theory are being introduced to the study of frames producing deep results in frame theory. At this very moment, the theory is beginning to grow rapidly with the host of new people entering the area. We will try to take a look at some of the current directions of research in abstract frame theory with an emphasis on the available techniques and unsolved problems. Some parts of the theory are so extensively developed, such as exponential (Fourier) frames, that we have chosen not to delve into them since even a small introduction which contains the essential tools and important results would be excessively long. This does not mean that these areas are not important, but to the contrary, they are so important as to require their own survey. Also, the bulk of the work on frame theory until recently has been in the applied directions, and we are not covering these areas at all (See the Remark below).

One of the nice things about frame theory is the fact that big portions are still underdeveloped - such as frames for finite dimensional Hilbert spaces, wavelet frames etc. Also, many of the extensively developed areas, such as Weyl-Heisenberg frames and exponential frames, still have many fundamental open questions to challenge anyone - such as the complete classification of Weyl-Heisenberg frames or the classification of exponential frames. Another interesting feature of the area is the broad spectrum of people working in different parts of it including Biologists, Engineers, Mathematicians (with backgrounds in Functional and Harmonic Analysis, Banach space theory, Operator Theory etc.), Physicists and more. Although each group has its own interests, there is an opportunity here to interact with a broad spectrum of researchers.

Now let us discuss the organization of the material. Section 2 contains all the background material needed throughout the paper, including results from Harmonic Analysis, Banach space theory, Operator Theory, and Hilbert space theory. Section 3 contains an (historical) introduction to frames through the work of Gabor and Duffin and Schaeffer. Section 4 is an introduction to abstract frame theory. This is one of the most extensive sections. It includes most of the fundamental results on abstract frame theory. Here we also make our case for treating frames as operators on a Hilbert space. That is, for bringing the power of operator theory into the frame theory arena. Section 5 is an introduction to Weyl-Heisenberg frames, one of the most important applications of frame theory. Here we look at the standard material on WH-frames including the important Zak Transform and its use to produce examples. This topic is taken up again in Section 9 where we look at some recent results in this area. Section 6 is a small introduction to perturbation theory containing some of the basic perturbation results. Section 7 is an introduction to wavelet frames. Section 8 concerns frames of translates and the related notions of exponential frames and sampling theory. Section 10 covers some selected topics in abstract frame theory including the projection methods, frames containing Riesz bases, and frames for Banach spaces.

There is a very well written introduction to WH-frames due to Heil and Walnut [96] which includes many of the results in this area through 1990. The author used this paper to enter the area of Weyl-Heisenberg frame theory and its influence can be seen throughout these notes. Also, the author wishes to thank Chris Heil for his careful reading of these notes and his many recommendations for imporvements. Let us also remark that until recently, many results in this area were "folklore" and therefore it is it is sometimes difficult or impossible to give proper credit for them.

**Remark:** These notes present a survey of modern abstract frame theory. However, until recently the majority of the work in frame theory was in the applied directions. Therefore, someone who only reads these notes will not be adequately introduced to the traditional major figures and traditional major topics in this area. We strongly recommend supplementing these notes with the excellent books by Daubechies [55], Feichtinger and Strohmer [74] and Young [139]. For a more accessible introduction to Weyl-Heisenberg frames, the forthcoming book of Gröchenig [88] is perfect. For a yet more applied approach to the area, we recommend Coifman and Zeevi [51] The above referenced material will give one a more balanced view of this important subject.

## 2. Background Material

We use $\mathbb{N}, \mathbb{Z}, \mathbb{R}, \mathbb{C}, \mathbb{Q}$ to denote the natural numbers, integers, real numbers, complex numbers, and rational numbers respectively. If we want to choose a number $a$ which may come from either $\mathbb{R}$ or $\mathbb{C}$, we just call $a$ **a scalar**. The torus group is denoted by $\mathbb{T}$ and is given by

$$\mathbb{T} = \{z \in \mathbb{C} : |z| = 1\}.$$

We identify the circle $\mathbb{T}$ with the interval $[0, 1)$ via the standard map $\zeta \to e^{2\pi i \zeta}$.

Sequences and series with undefined limits are taken to be over $\mathbb{Z}$ and integrals with undefined limits are taken to be over $\mathbb{R}$. Integration is always with respect to Lebesgue measure. **All functions will be assumed to be measurable**. A property is said to hold **almost everywhere**, denoted a.e., if the set of points in $\mathbb{R}$ (or $\mathbb{C}$) where it does not hold has Lebesgue measure zero. All functions (unless otherwise indicated) have domain the real line and take values in $\mathbb{C}$. The **characteristic function** of a set $E \subset \mathbb{R}$ is

$$\chi_E(x) = \begin{cases} 1 : & x \in E \\ 0 : & x \notin E \end{cases}$$

The **Kronecker delta** defined for $m, n \in \mathbb{Z}$ is

$$\delta_{mn} = \begin{cases} 1: & m = n \\ 0: & m \neq n \end{cases}$$

The **essential supremum** of a function $f$ is

$$\|f\|_\infty = \text{ess sup}_{x \in \mathbb{R}} |f(x)| = \inf\{\lambda \in \mathbb{R} : f(x) \leq \lambda \ a.e.\}.$$

We will work in general with the Banach spaces $L^p(\mathbb{R})$, $1 \leq p \leq \infty$. For $p = \infty$, this is the space (with norm $\|\cdot\|_\infty$)

$$L^\infty(\mathbb{R}) = \{f : \mathbb{R} \to \mathbb{C} : \|f\|_\infty < \infty\}.$$

For $1 \leq p < \infty$ we have the Banach space (with norm $\|f\|_p$)

$$L^p(\mathbb{R}) = \{f : \mathbb{R} \to \mathbb{C} : \|f\|_p = \left(\int_{\mathbb{R}} |f(x)|^p dx\right)^{1/p} < \infty\}.$$

For any set of vectors $E$ in a linear space $X$, we write span E to denote the family of finite linear combinations of elements of $E$. Or if E is a sequence of elements, say $(x_n)$, we write span E $= \text{span}_n(x_n)$. The closure of the span of $E$ will be denoted $\overline{span} \ E$. If the span of $(x_n)$ is dense in $X$ we say that $(x_n)$ is **complete**. We will also need the atomic version of $L^2(\mathbb{R})$ denoted $\ell_2$. This is the space of all sequences of scalars $(a_n)$ for which

$$\|(a_n)\|_{\ell_2} = \left(\sum_n |a_n|^2\right)^{1/2} < \infty.$$

Since we are working only in separable spaces, we can consider all infinite dimensional **Hilbert spaces** as $L^2(\mathbb{R})$ or $\ell_2$ and will denote it abstractly as the **Hilbert space H**. For a natural number n, we denote an n-dimensional Hilbert space by $H_n$. We assume the reader is familiar with the basic properties of a Hilbert space. We will denote the **inner product** for $H$ by $< \cdot, \cdot >$ and recall that $\|x\|^2 = <x, x>$, for all $x \in H$. Given $E \subset H$ and $F \subset H$, we say that $E$ is **orthogonal** to $F$ (and write $E \perp F$) if $<x, y> = 0$, for all $x \in E, y \in F$. If $E \subset H$, we denote the **orthogonal complement** of $E$ by

$$E^\perp = \{x \in H :<x, y> = 0, \text{ for all } y \in E\}.$$

If $(x_n)$ is a sequence of vectors in $H$, we say $(x_n)$ is an **orthogonal sequence** if $<x_n, x_m> = 0$, for all $n \neq m$. If moreover, $\|x_n\|^2 = <x_n, x_n> = 1$, we call $(x_n)$ an **orthonormal sequence**. We generally denote an orthonormal sequence by $(e_n)$. If furthermore the sequence is complete, we call $(x_n)$ an **orthonormal basis** for $H$. A pair of sequences $(x_i), (y_i)$ in $H$ is called a **biorthogonal system** if $<x_i, y_j> = 0$, for all $i \neq j$.

For an orthonormal sequence $(x_n)$ in a Hilbert space $H$, the following are equivalent:

(1) $(e_n)$ is complete.

(2) We have the **Plancherel formula**:

$$\|x\|^2 = \sum_n | < x, e_n > |^2,$$

for all $x \in H$.

(3) $x = \sum_n < x, e_n > e_n$, for all $x \in H$.

In this case we call $(< x, e_n >)$ the **Fourier coefficients** of $x$ (with respect to the orthonormal basis $(e_n)$). It follows that the coefficients $(< x, e_n >)$ in (3) are unique. This is in contrast to a frame which we will see may have infinitely many different representations in general.

Let $H, K$ be Hilbert spaces with inner products $< \cdot, \cdot >_H, < \cdot, \cdot >_K$ and norms $\| \cdot \|_H, \| \cdot \|_K$ respectively and let $T : H \to K$. We define the following properties.

(1) $T$ is **linear** if $T(ax + by) = aTx + bTy$, for all scalars $a, b$ and all $x, y \in H$.

(2) $T$ is 1-1 (or **injective**) if $Tx \neq Ty$, for all $x \neq y$.

(3) The **range** of $T$ is Rng $T = \{Tx; x \in H\}$ and the **rank** of $T$ is the dimension of Rng T. The **co-rank** of $T$ is the dimension of $(\text{Rng } T)^{\perp}$.

(4) The **kernel** of $T$ is ker $T = \{x : Tx = 0\}$ and the **nullity** of $T$ is the dimension of ker T.

(5) $T$ is **onto** (or **surjective**) if Rng T = K.

(6) The norm of $T$ is

$$\|T\| = \sup_{0 \neq x \in H} \frac{\|Tx\|}{\|x\|} = \sup_{\|x\|_H = 1} \|Tx\|.$$

We say that $T$ is **bounded** if $\|T\| < \infty$. A linear operator $T$ is bounded if and only if it is continuous.

From here on we will assume that $T$ is a bounded linear operator.

(7) The **adjoint** of $T$ is the unique operator $T^* : K \to H$ satisfying

$$< Tx, y >_K = < x, T^*y >_H, \text{ for all } x \in H, y \in K.$$

It is a simple calculation to show that $\|T\| = \|T^*\|$.

(8) We say that $T$ is an **isomorphism** if it is 1-1, continuous and it has an inverse $T^{-1}$ (defined on Rng $T$) which is continuous. This is equivalent to the existence of a constant $M > 0$ satisfying:

$$\frac{1}{M}\|x\| \leq \|Tx\| \leq M\|x\|, \text{ for all } x \in H.$$

If $T$ is also onto, we say that $T$ is an **invertible operator**. If $T : H \rightarrow H$ satisfies $\|I - T\| < 1$, then $T$ is invertible and its inverse can be represented by the **Neuman series** [119]

$$T^{-1}f = \sum_{n=0}^{\infty}(I - T)^n f.$$

(9) $T$ is an **isometry** if $\|Tx\|_K = \|x\|_H$, for all $x \in H$. It can be shown that $T$ is an isometry if and only if $< Tx, Ty >_K = < x, y >_H$, for all $x, y \in H$. $T$ is a **co-isometry** if its adjoint is an isometry.

(10) $T$ is a **unitary map** if it is an invertible isometry.

(11) $T$ is a **partial isometry** if $T$ is an isometry on the orthogonal complement of its kernel.

From here on, we will let $H = K$ and let $T, S : H \rightarrow H$ be bounded linear operators and let $(e_n)$ be an orthonormal basis for $H$.

(12) $T$ is **self-adjoint** if $T = T^*$. This is equivalent to:

$$< Tx, y >=< x, Ty >, \text{ for all } x, y \in H.$$

(13) $T$ is **positive**, denoted $T \geq 0$, if $< Tx, x > \geq 0$, for all $x \in H$. Positive operators are self-adjoint.

(14) We write $T \geq S$ if $T - S \geq 0$. If $T$ is positive and $T \leq I$ then we can write

$$T = \frac{1}{2}(W + W^*),$$

where $W$ is a unitary operator.

(15) There exists a partial isometry $U$ on $H$ and a positive operator $P$ on $H$ such that $T = UP$. We can find the operators $U, P$ so that $\ker U = \ker P$, and this additional condition uniquely determines them. This representation of $T$ is called the **polar decomposition** of $T$. A necessary and sufficient condition that $U$ be an isometry is that $T$ be 1-1, and a necessary and sufficient condition that $U$ be a co-isometry is that Rng $T$ be dense in $H$.

(16) The **trace** of $T$ is given by: $\text{Tr } T = \sum_n < Te_n, e_n >$.

(17) The **Hilbert-Schmidt norm** of $T$ is

$$\|T\|_{HS} = \left( \sum_n \|Te_n\|^2 \right)^{1/2}.$$

If $P$ is a projection of rank n, then

$$\|P\|_{HS}^2 = \sum_i \|Pe_i\|^2 = \sum_i < Pe_i, Pe_i >$$

$$= \sum_i < Pe_i, e_i >= \text{Tr } T = n.$$

For a function $f$ on $\mathbb{R}$ we define the operators:

$$\begin{aligned}
\text{Translation:} \quad & T_a f(x) = f(x-a), & a \in \mathbb{R} \\
\text{Modulation:} \quad & E_a f(x) = e^{2\pi i a x} f(x), & a \in \mathbb{R} \\
\text{Dilation:} \quad & D_a f(x) = |a|^{-1/2} f(x/a), & a \in \mathbb{R} - \{0\}
\end{aligned}$$

We also use the symbol $E_a$ to denote the **exponential function** $E_a(x) = e^{2\pi i a x}$. The **two-dimensional exponentials** are $E_{a,b}(x,y) = e^{2\pi i a x} e^{2\pi i b y}$. Each of the operators $T_a, E_a, D_a$ are unitary operators on $L^2(\mathbb{R})$ and they satisfy:

$$\begin{aligned}
T_a E_b f(x) &= e^{2\pi i b(x-a)} f(x-a); \\
E_b T_a f(x) &= e^{2\pi i b x} f(x-a); \\
T_b D_a f(x) &= |a|^{-1/2} f(\tfrac{x-b}{a}); \\
D_a T_b f(x) &= |a|^{-1/2} f(\tfrac{x}{a} - b); \\
E_b D_a f(x) &= e^{2\pi i b x} |a|^{-1/2} f(\tfrac{x}{a}); \\
D_a E_b f(x) &= e^{2\pi i b x/a} |a|^{-1/2} f(\tfrac{x}{a}).
\end{aligned}$$

The **Fourier transform** of a function $f \in L^1(\mathbb{R})$ is

$$\hat{f}(\gamma) = \int_R f(x) e^{-2\pi i \gamma x} dx, \text{ for all } \gamma \in \mathbb{R}.$$

We define the Fourier transform of functions $f \in L^2(\mathbb{R})$ by observing that the above definition works on a dense subspace of $L^2(\mathbb{R})$ and continuing it to the closure. We have

$$\widehat{T_a f} = E_{-a}\hat{f}; \quad \widehat{E_a f} = T_a \hat{f}; \quad \widehat{D_a f} = D_{1/a}\hat{f}$$

We also have the **Plancherel formula**

$$\|f\|_2 = \|\hat{f}\|_2, \text{ for all } f \in H,$$

and the **Parseval formula**

$$< f, g > = < \hat{f}, \hat{g} > . \text{ for all } f, g \in H.$$

We end with a discussion of the types of convergence we will be working with. These results can be found in [62, 109]. Given two sequences $(x_i), (y_i)$ spanning (Banach) Hilbert spaces $X, Y$, respectively, we say the sequences are **equivalent** and write $(x_i) \approx (y_i)$ if setting $Tx_i = y_i$ uniquely extends to a well-defined invertible operator from $X$ onto $Y$. We say that $(x_i)$ is a **Schauder basis** (or just a **basis**) for $X$ if every element $x \in X$ has a unique representation in the form,

$$(2.1) \qquad x = \sum_i a_i x_i,$$

where $(a_i)$ is a sequence of scalars. Associated to a basis $(x_i)_{i \in \mathbb{N}}$ is the **basis projections** $(P_n)$ given by:

$$P_n(\sum_{i=1}^{\infty} a_i x_i) = \sum_{i=1}^{n} a_i x_i.$$

We define the **basis constant** of $(x_i)$ to be $\sup_n \|P_n\|$. If the series in (2.1) converge unconditionally (see below) for all $x \in X$ we call $(x_i)$ an **unconditional basis** for $X$. There is also a **unconditional basis constant** which we will not define in general. The basis $(x_i)$ is **bounded** if there is a constant $M > 0$ so that $M^{-1} \le \|x_i\| \le M$, for all $i$. A **Riesz basis** for a Hilbert space $H$ is a bounded, unconditional basis for $H$. It is known (see [109]) that $(x_i)$ is a Riesz basis for $H$ if and only if $(x_i) \approx (e_i)$, where $(e_i)$ is an orthonormal basis for $H$. In this case, we define the **Riesz basis constants** for $(x_i)$ to be the largest number $K$ and the smallest number $M$ satisfying, for all sequences of scalars $(a_i)$,

$$K \left( \sum_i |a_i|^2 \right)^{1/2} \le \| \sum_i a_i x_i \| \le M \left( \sum_i |a_i|^2 \right)^{1/2}.$$

We will be working with several different forms of convergence for series in a (Banach) Hilbert space $X$. We say that a series $\sum_i x_i$ of elements of $X$ **converges unconditionally** if every rearrangement of $(x_i)$ converges to an element of $X$ (and therefore to the same element of $X$). It is known (see for example [54, 109])

**Theorem 2.1.** *For $x_i$ in a Banach space $X$, the following are equivalent:*
(1) $\sum_i x_i$ *converges unconditionally in $X$.*

(2) *For every increasing sequence of natural numbers $(n_i)$ we have that $\sum_i x_{n_i}$ converges in $X$.*

(3) $\sum_i \theta_i x_i$ *converges for every choice of $\theta_i = \pm 1$.*

*Moreover, in this case, there is a constant $M > 0$ so that for all choices of scalars $(a_i)$ we have*

$$\|\sum_i a_i x_i\| \le M \sup_i |a_i| \cdot \|\sum_i x_i\|.$$

We will also work with **weak convergence** in our spaces. In a Hilbert space we say that $(x_i)$ **converges weakly** to $x \in H$ if

$$\lim_i < y, x_i > = < y, x >, \text{ for all } y \in H.$$

Again, if for every rearrangement $(x_{\sigma(i)})$ of $(x_i)$ the series $\sum_i x_{\sigma(i)}$ converges weakly to an element of $H$, we say that $\sum_i x_i$ is **weakly unconditionally convergent**. The celebrated Orlicz-Pettis Theorem says that weak unconditional convergence is the same as unconditional convergence in every Banach space.

**Theorem 2.2. (Orlicz-Pettis Theorem)** *In a (Banach) Hilbert space $H$, a series is weakly unconditionally convergent if and only if it is norm unconditionally convergent.*

Sometimes it is useful to know that a series is convergent without knowing what it is converging to. This is the notion of "Cauchy sequences" in $\mathbb{R}$. In a Banach space, this becomes **wuC**. A series $\sum x_i$ in a (Banach) Hilbert space $H$ is said to be **weakly unconditionally Cauchy** if given any permutation $\sigma$ of $\mathbb{N}$, $\sum_{i=1}^n x_{\sigma(i)}$ is a weakly Cauchy sequence in $H$. Recall the Banach space $c_0$:

$$c_0 = \{(a_i) : \|(a_i)\|_{c_0} = \sup_i |a_i| < \infty, \text{ and } \lim_i a_i = 0\}.$$

The following theorem clearly holds for a Hilbert space.

**Theorem 2.3.** *If $c_0$ does not embed into a Banach space $X$, then every series $\sum_i x_i$ which is wuC is also unconditionally convergent in $X$.*

## 3. Historical Development of Frames

To give us a frame of reference, we will start by attending a concert. Suppose, for a moment, you are sitting in the audience at a piano concert. You are enjoying the sound of the piano filling up the concert hall. What you are hearing is a continuous acoustical signal (if we ignore some technicalities such as the "striking" of the keys). But to the pianist, this is something else altogether. To this person, "the concert" is a book of sheet music which consists of sheets of paper with special lines and artistically placed black dots representing notes. The sheet music is a "discretization" of the musical piece, localized in time and frequency. The vertical direction represents frequencies and the horizontal direction represents time broken into equally spaced intervals. The sheet music tells the pianist which notes to play in which time intervals. Technically, if you were musically gifted, you could sit in the audience and write the sheet music for the piece you were hearing. Then you could bring this home and play the piece back (called "reconstruction of the signal"). More importantly, if you found something in the music irritating (like "noise" to your ear) you could erase those dots from the sheet music then play (reconstruct) a variation of the music with the "noise" removed. Although it would be difficult to read and play, technically we could write our sheet music by placing a number $\{0, 1, 2, 3\}$ in each note position vertically, equally spaced in time. These values would represent the "intensity" of that note as being: (0) don't play it, (1) play it, (2) play it flat, or (3) play it sharp respectively. We can interpret all of our work on (at least Weyl-Heisenberg) frames in the context of writing sheet music. Our sheet music will be a little more sophisticated in that time will continue forever through the past and the future. Also, we will have infinitely many choices of frequencies (which now will be called "elementary signals"), and for each frequency we may have any complex number for its intensity.

In 1946, D. Gabor [84] formulated a fundamental approach for signal decomposition in terms of elementary signals. Gabor's approach quickly became a paradigm for the spectral analysis associated with time-frequency methods, such as the short-time Fourier transform and the Wigner transform. Gabor's idea required a tiling of the time-frequency domain (also called the information plane or phase space depending upon the area you work in) by non-overlapping half open rectangles. Gabor reasoned that certain "optimal" elementary signals should provide an efficient decomposition of the information plane. This decomposition should determine countably many components of the information plane with each component sufficiently localized in time and frequency so that a coefficient c associated with a component R would characterize the amount of information from R in the signal. Moreover, it should not require smaller components to distinguish different types of information from R.

So Gabor needed to find a countable collection of "optimal" elementary signals with small associated areas in the information plane. Gabor chose modulates and translates of Gaussians as elementary signals because the product of their time and frequency variances are optimal relative to the Classical Uncertainty Principle Inequality (See Section 5).

So how does Gabor's scheme work? Gabor let

$$g(t) = \pi^{-1/4} e^{-t^2/2} \in L^2(\mathbb{R}).$$

The function $g$ is called the **window function**. We fix $a, b \in \mathbb{R}^+$. Our elementary signals are then

$$(E_{mb}T_{na}g)_{m,n \in \mathbb{Z}}.$$

We choose a signal $f$ (i.e. a function $f \in L^2(\mathbb{R})$) and compute the **intensity** of each elementary signal in $f$. To do this we fix $n = 0$ and compute the inner products

$$< f, E_{mb}g > = < f \cdot g, E_{mb} >, \quad \text{for all} \quad m \in \mathbb{Z}.$$

So we are just computing Fourier coefficients for $f \cdot g$ (or the "weighted" Fourier transform of f). These coefficients give an indication of the frequency content of the signal $f$ in a neighborhood of 0. Now, translate the window by say $n = 1$ and do it all over again,

$$< f, E_{mb}T_a g > = < f \cdot T_a g, E_{mb} >, \quad \text{for all} \quad m \in \mathbb{Z}.$$

We continue to compute these Fourier coefficients for all $n \in \mathbb{Z}$. Now we have a set of coefficients $(c_{mn})_{m,n \in \mathbb{Z}}$ associated with the elementary signals which should be unique to our signal (the exact requirements on $g, a, b$ for uniqueness is a deep question which is covered in detail in Sections 5, 9). To return to our opening example, the elementary signals are our "notes" and the $c_{mn}$'s are the intensities of these notes in our signal. This association is unique and allows us to work with the signal as just a discrete set of intensities associated with the elementary signals.

We can view Gabor's method as a sampled short time Fourier transform (i.e. a set of inner products of f with elementary signals $(E_m T_n g)$ with discrete labels in the lattice $a\mathbb{Z} \times b\mathbb{Z}$). Recall that the **short time Fourier transform** of $f \in L^2(\mathbb{R})$ (with respect to a window function $g$) is

$$\nu_g(t, \gamma) = \int_{\mathbb{R}} f(s)\overline{g(t-s)}e^{-2\pi i \gamma s} \ ds = < f, T_t M_\gamma g > .$$

Also, we can recover $f$ from its short time Fourier transform by the inversion formula

$$f(t) = \frac{1}{\|g\|_{L^2}^2} \int_{\mathbb{R} \times \mathbb{R}} \nu_g f(s, \gamma)g(t-s) \ e^{2\pi i \gamma t} \ dt \ d\gamma.$$

In light of this, the sampled short time Fourier transform is also referred to as the **Gabor transform**.

Gabor's use of Gaussians has both advantages and disadvantages which will be discussed in detail later. We will not go any further with this example at this time since we will cover these frames in full generality in Section's 5 and 9.

In 1952, Duffin and Schaeffer [63] were working on some deep problems in nonharmonic Fourier series. Much of the work in nonharmonic Fourier series was initiated by the fundamental results of Paley and Wiener [121]. Duffin and Schaeffer were working with families of exponentials $(e^{i\lambda_n t})_{n \in \mathbb{Z}}$ trying to determine when they are complete or form a Riesz basis for $L^2[a, b]$ etc. This led them to define

**Definition 3.1.** *A sequence* $(f_n)_{n \in \mathbb{Z}}$ *of elements of a Hilbert space* $H$ *is called a* **frame** *if there are constants* $A, B > 0$ *such that*

$$(3.1) \qquad A\|f\|^2 \leq \sum_{n \in \mathbb{Z}} | < f, f_n > |^2 \leq B\|f\|^2, \quad \text{for all} \ \ f \in H.$$

One of Duffin and Schaeffer's main results (see [63], Theorem 1) is that $(e^{-2\pi i \lambda_n t})_{n \in \mathbb{Z}}$ is a frame for $L^2[-\gamma, \gamma]$ if $(\lambda_n)_{n \in \mathbb{Z}}$ is uniformly dense with uniform density greater than $2\gamma$. Uniform density and the role played in irregular sampling is discussed in [15]. These frames are called **exponential frames** or **Fourier frames** and will be discussed again in Section 8

For some reason, the notion of a frame introduced by Duffin and Schaeffer was not followed up on outisde of nonharmonic Fourier series. However, it was brought back to life in 1986 by Daubechies, Grossman and Meyer [57] right at the dawn of the Wavelet era.

## 4. An Introduction to Frame Theory

Good references for this section are [93] and [139]. The numbers $A, B$ in the definition of a frame, Definition 3.1, are called the **lower** and **upper frame bounds** respectively. The largest number $A > 0$ and smallest number $B > 0$ satisfying the frame inequalities for all $f \in H$ are called the **optimal frame bounds**. The frame is a **tight frame** if $A = B$ and a **normalized tight frame** if $A = B = 1$. A frame is **exact** if it ceases to be a frame when any one of its elements is removed. As we will see, a frame is exact if and only if it is a Riesz basis. A non-exact frame is called **overcomplete** in the sense that at least one vector can be removed from the frame and the remaining set of vectors will still form a frame for $H$ (but perhaps with different frame bounds).

It is immediate that an orthonormal basis $(e_n)$ for $H$ is a normalized tight frame for $H$. But the following sequences are also normalized tight frames for $H$:

$$\{e_1, 0, e_2, 0, e_3, 0 \cdots\},$$

$$\{\frac{e_1}{\sqrt{2}}, \frac{e_1}{\sqrt{2}}, \frac{e_2}{\sqrt{2}}, \frac{e_2}{\sqrt{2}}, \cdots\},$$

$$\{e_1, \frac{e_2}{\sqrt{2}}, \frac{e_2}{\sqrt{2}}, \frac{e_3}{\sqrt{3}}, \frac{e_3}{\sqrt{3}}, \frac{e_3}{\sqrt{3}} \cdots\}.$$

If $(f_n)$ is any sequence with a finite upper frame bound, then $(e_n, f_n)$ is a frame for $H$. It is obvious that $(e_n/n)$ fails to have a lower frame bound while $(ne_n)$ fails to have a finite upper frame bound.

For the $n$-dimensional Hilbert space $H_n$ with orthonormal basis $(e_i^n)_{i=1}^n$, a direct calculation shows that the following sequence is a normalized tight frame for $H_n$:

$$(4.1) \qquad f_j^n = e_j^n - \frac{1}{n} \sum_{i=1}^n e_i^n \ \ \text{for} \ \ j = 1, 2, \cdots, n,$$

$$f_{n+1}^n = \frac{1}{\sqrt{n}} \sum_{i=1}^n e_i^n.$$

We can also see that $(f_j^n)_{j=1}^{n+1}$ is a normalized tight frame for $H_n$ by applying Theorem 4.10 below. That is, if $P$ is the rank one orthogonal projection on $H_n$ given by:

$$P(\sum_{i=1}^n a_i e_i^n) = (\sum_{i=1}^n \frac{a_i}{n}) \sum_{i=1}^n e_i^n,$$

then $P^\perp e_j^n = f_j^n$.

One fruitful approach to frame theory for infinite dimensional Hilbert spaces is to view frames as operators. We feel that this is the best approach since it brings the tremendous power of operator theory, $C^*$-algebras etc. to bear on the subject. A good introduction to this approach can be found in the Memoir of Han and Larson [93]. To formulate this approach, let $(e_n)$ be an orthonormal basis for an infinite dimensional Hilbert space $H$ and let $f_n \in H$, for all $n \in \mathbb{Z}$. We call the operator $T : H \to H$ given by $Te_n = f_n$ the **preframe operator** associated with $(f_n)$. Now, for each $f \in H$ and $n \in \mathbb{Z}$ we have $< T^*f, e_n > = < f, Te_n > = < f, f_n >$. Thus

$$(4.2) \qquad T^*f = \sum_n < f, f_n > e_n, \ \ \text{for all} \ \ f \in H.$$

It follows that the preframe operator is bounded if and only if $(f_n)$ has a finite upper frame bound $B$. Also, by (4.2)

$$\|T^*f\|^2 = \sum_n | < f, f_n > |^2, \ \ \text{for all} \ \ f \in H.$$

Comparing this to Definition 3.1 we have

**Theorem 4.1.** *Let $H$ be a Hilbert space with an orthonormal basis $(e_n)$. Also let $(f_n)$ be a sequence of elements of $H$ and let $Te_n = f_n$ be the preframe operator. The following are equivalent:*

*(1) $(f_n)$ is a frame for $H$.*

*(2) The operator $T$ is bounded, linear and onto.*

*(3) The operator $T^*$ is an (possibly into) isomorphism (called the **frame transform**).*

*Moreover, $(f_n)$ is a normalized tight frame if and only if the preframe operator is a quotient map (i.e. a partial isometry).*

The dimension of the kernel of T is called the **excess** of the frame. It follows that $S = TT^*$ is an invertible operator on $H$, called the **frame operator**. Moreover, we have

$$Sf = TT^*f = T(\sum_n <f, f_n> e_n) = \sum_n <f, f_n> Te_n = \sum_n <f, f_n> f_n.$$

A direct calculation now yields

$$<Sf, f> = \sum_n | <f, f_n> |^2.$$

Therefore, the **frame operator is a positive, self-adjoint invertible operator** on $H$. Also, the frame inequalities (3.1) yield that $(f_n)$ is a frame with frame bounds $A, B > 0$ if and only if $A \cdot I \leq S \leq B \cdot I$. Hence, $(f_n)$ is a normalized tight frame if and only if $S = I$. Also, a direct calculation yields

$$(4.3) \qquad f = SS^{-1}f = \sum_n <S^{-1}f, f_n> f_n$$

$$= \sum_n <f, S^{-1}f_n> f_n = \sum_n <f, S^{-1/2}f_n> S^{-1/2}f_n.$$

We call $(<S^{-1}f, f_n>)$ the **frame coefficients** for $f$. It follows

**Theorem 4.2.** *Every frame $(f_n)$ (with frame operator $S$) is equivalent to the normalized tight frame $(S^{-1/2}f_n)$.*

Note that equation (4.3) is our "reconstruction formula" for an element $f \in H$. This points out one of the major problems encountered in applications of frame theory. In order to reconstruct a vector we need to find its frame coefficients, i.e. we first have to find $S^{-1}f$. This requires inverting an infinite matrix - no simple task. Since $S$ is an isomorphism on $H$, $(S^{-1}f_n)$ is a frame equivalent to the frame $(f_n)$ and is called **the (canonical) dual frame**.

**Remark:** This approach to frame theory for infinite dimensional spaces also works, with a slight variation, for finite dimensional spaces. A sequence $(f_n)$ in a m-dimensional Hilbert space $H_m$ is a frame if and only if the operator $T : \ell_2 \to H_m$ given by $Te_n = f_n$ is bounded, linear and onto, where $(e_n)$ is any orthonormal basis for $\ell_2$. We will almost exclusively work in infinite-dimensional separable Hilbert spaces except in a few cases where we will make it clear that we are working in $H_m$.

Now we will consider some of the basic properties of frames.

**Proposition 4.3.** *Let $(f_n)$ be a frame for a Hilbert space $H$. Then the following are equivalent:*
*(1) $(f_n)$ is exact.*

*(2) $(f_n)$ is a Riesz basis.*

*(3) $(f_n)$ is a frame which is $\omega$-independent.*

*Moreover, in this case, the Riesz basis constants for $(f_n)$ are the numbers $\sqrt{A}$, $\sqrt{B}$ where $A, B$ are the frame bounds.*

*Proof.* $(1) \Leftrightarrow (2)$ We note that $(f_n)$ is exact if and only if the preframe operator $T$ is one-to-one. But since $T$ is bounded, linear and onto, this happens if and only if $T$ is an invertible operator - which makes $(f_n)$ a Riesz basis.

$(1) \Leftrightarrow (3)$ This is similar.

We leave the moreover part of the theorem to the reader. $\qquad \square$

Kim and Lim [106] have given equivalent conditions for a frame to be a Riesz basis and compute the bounds in terms of the eigenvalues of the Gram matrices of finite subsets.

There is a simple method for showing that two frames are equivalent.

**Proposition 4.4.** *Let $(f_i)$ and $(g_i)$ be frames for a Hilbert space $H$ with pre-frame operators $T_1, T_2$ respectively with respect to a fixed orthonormal basis $(e_i)$. For each finitely non-zero sequence of scalars $(a_i)$ let $L(\sum_i a_i f_i) = \sum_i a_i g_i$ be a relation. The following are equivalent:*

*(1) $L$ is a function.*

*(2) $L$ is a bounded linear operator.*

*(3) $\ker T_1 \subset \ker T_2$.*

*Proof.* $(1) \Rightarrow (3)$: Since $L$ is linear, $L$ is a function implies $L(0) = 0$.

$(3) \Rightarrow (2)$: We get immediately from $(3)$ that $L$ is a function which is clearly linear. To see that $L$ is bounded we note that

$$\|L \sum_i a_i f_i\| = \|\sum_i a_i g_i\| = \|T_2(\sum_i a_i e_i)\| \approx \|P_{(ker\ T_2)^\perp}(\sum_i a_i e_i)\|$$

$$\leq \|P_{(ker\ T_1)^\perp} \sum_i a_i e_i\| \approx \|T_1(\sum_i a_i e_i)\| = \|\sum_i a_i f_i\|.$$

$(2) \Rightarrow (1)$: This is obvious. $\qquad \square$

**Corollary 4.5.** *If $(f_i)$ and $(g_i)$ are frames for a Hilbert space $H$ with preframe operators $T_1, T_2$ respectively, the following are equivalent:*

*(1) $(f_i) \approx (g_i)$.*

*(2) $\ker T_1 = \ker T_2$.*

*(3) For all sequences of scalars $(a_i)$ we have that $\sum_i a_i f_i = 0$ if and only if $\sum_i a_i g_i = 0$.*

Our next proposition shows the relationship between the frame elements and the frame bounds.

**Proposition 4.6.** *Let $(f_n)$ be a frame for $H$ with frame bounds $A, B$. We have for all $n \in \mathbb{Z}$ that $\|f_n\|^2 \leq B$ and $\|f_n\|^2 = B$ implies $f_n \perp span_{j\neq n}f_j$. If $\|f_n\|^2 < A$, then $f_n \in \overline{span}(f_j)_{j\neq n}$.*

*Proof.* If we replace $f$ in Definition 3.1 by $f_n$ we see that

$$A\|f_n\|^2 \leq \|f_n\|^4 + \sum_{j\neq n} | < f_n, f_j > |^2 \leq B\|f_n\|^2.$$

The first part of the result is now immediate. For the second part, assume to the contrary that $E = \text{span}(f_j)_{j\neq n}$ is a proper subspace of $H$. Replacing $f_n$ in the above inequality by $P_{E^\perp}f_n$ and using the left hand side of the inequality yields an immediate contradiction. ∎

As a particular case of Proposition 4.6 we have for a normalized tight frame $(f_n)$ that $\|f_n\|^2 \leq 1$ and $\|f_n\| = 1$ if and only if $f_n \perp span_{j\neq n}f_j$. If $(f_n)$ is an exact frame then $< S^{-1}f_n, f_m > = < S^{-1/2}f_n, S^{-1/2}f_m > = \delta_{nm}$ (where $\delta_{nm}$ is the Kronecker delta) since $(S^{-1/2}f_n)$ is now an orthonormal basis for $H$. That is, $(S^{-1}f_n)$ and $(f_n)$ form a biorthogonal system. Also, it follows that $(e_n)$ is an orthonormal basis for $H$ if and only if it is an exact, normalized tight frame. Another consequence of Proposition 4.6 is

**Proposition 4.7.** *The removal of a vector from a frame leaves either a frame or an incomplete set.*

*Proof.* By Theorem 4.2, we may assume that $(f_i)$ is a normalized tight frame. Now, by Proposition 4.6, for any n, either $\|f_n\| = 1$ and $f_n \perp span_{j\neq n}f_j$, or $\|f_n\| < 1$ and $f_n \in span_{j\neq n}f_j$. ∎

Since a frame is not $\omega$-independent (unless it is a Riesz basis) a vector in the space may have many representations relative to the frame besides the natural one given by the frame coefficients. However, the natural representation of a vector is the unique representation of minimal $\ell_2$-norm as the following result of Duffin and Schaeffer [63] shows.

**Theorem 4.8.** *Let $(f_n)$ be a frame for a Hilbert space $H$ and $f \in H$. If $(b_n)$ is any sequence of scalars such that*

$$f = \sum_n b_n f_n,$$

*then*

$$(4.4) \qquad \sum_n |b_n|^2 = \sum_n | < S^{-1}f, f_n > |^2 + \sum_n | < S^{-1}f, f_n > -b_n|^2.$$

*Proof.* We have by assumption

$$\sum_n < S^{-1}f, f_n > f_n = \sum_n b_n f_n.$$

Now, taking the inner product of both sides with $S^{-1}f$ we get

$$\sum_n |<S^{-1}f, f_n>|^2 = \sum_n \overline{<S^{-1}f, f_n>} b_n,$$

and (4.4) follows easily. $\qquad\square$

A major advantage of frames over wavelets (see Section 7) is that orthogonal projections take frames to frames but do not map wavelets to wavelets.

**Proposition 4.9.** *Let $(f_n)$ be a frame for $H$ with frame bounds $A, B$, and let $P$ be an orthogonal projection on $H$. Then $(Pf_n)$ is a frame for $P(H)$ with frame bounds $A, B$.*

*Proof.* For any $f \in P(H)$ we have

$$\sum_n |<f, Pf_n>|^2 = \sum_n |<Pf, f_n>|^2 = \sum_n |<f, f_n>|^2.$$

The result is now immediate. $\qquad\square$

It follows that an orthogonal projection $P$ applied to an orthonormal basis $(e_n)$ (or just a normalized tight frame) yields a normalized tight frame $(Pe_n)$ for $P(H)$. The converse of this is also true and is a result of Han and Larson [93].

**Theorem 4.10.** *A sequence $(f_n)$ is a normalized tight frame for a Hilbert space $H$ if and only if there is a larger Hilbert space $H \subset K$ and an orthonormal basis $(e_n)$ for $K$ so that the orthogonal projection $P_H$ of $K$ onto $H$ satisfies $Pe_n = f_n$ for all $n = 1, 2, \cdots$.*

*Proof.* The "only if" part follows from Proposition 4.9. For the "if" part, if $(f_n)$ is a normalized tight frame for $H$ then the preframe operator $T : \ell_2 \to H$ is a partial isometry. Let $(e_n)$ be an orthonormal basis for $\ell_2$ for which $T(e_n) = f_n$ is our frame. Since $T^*$ is an into isometry we can associate $H$ with $T^*(H)$. Now let $K = \ell_2$ and $P$ be the orthogonal projection of $K$ onto $T^*(H)$. Then for all $n = 1, 2, \cdots$ and all $g = T^*f \in T^*(H)$ we have

$$<T^*f, Pe_n> = <T^*f, e_n> = <f, Te_n> = <f, f_n> = <T^*f, T^*f_n>.$$

It follows that $Pe_n = T^*f_n$, and by our association of $H$ with $T^*(H)$, $(T^*f_n)$ is our frame. $\qquad\square$

Casazza, Han and Larson [41] have generalized Theorem 4.10 to show that any sequence $(f_n)$ in a Hilbert space $H$ is a frame for $H$ if and only if there is a larger Hilbert space $H \subset K$, an orthonormal basis $(e_n)$ for $K$ and a (not necessarily orthogonal) projection $P : K \to H$ for which $Pe_n = f_n$, for all $i$.

Another interesting consequence of considering frames as operators, in the finite-dimensional case, is the following finite-dimensional result.

**Proposition 4.11.** *If $(f_n)$ is a normalized tight frame for an $m$-dimensional Hilbert space $H_m$ then*

$$\sum_n \|f_n\|^2 = m.$$

*Hence, if $(f_n)$ is a normalized tight frame for an infinite dimensional Hilbert space $H$, then for all finite rank orthogonal projections $P$ on $H$ we have that*

$$\sum_n \|P f_n\|^2 \quad \text{is a natural number}$$

*Proof.* By Theorem 4.10, there is a larger Hilbert space $H \subset K$ and an orthogonal projection $P$ taking an orthonormal basis $(e_n)$ for $H$ to $P e_n = f_n$. Now,

$$\sum_n \|f_n\|^2 = \sum_n \|P e_n\|^2 = \|P\|_{HS}^2 = \dim H_m,$$

where $\| \cdot \|_{HS}$ denotes the Hilbert Schmidt norm of $P$.                    $\square$

Proposition 4.11 can be generalized to an arbitrary frame. Recall that a subspace of codimension 1 in a Hilbert space is called a **hyperplane**.

**Proposition 4.12.** *Let $(f_i)$ be a sequence in a Hilbert space $H$. The optimal lower (respectively, upper) frame bound $A$ (resectively, $B$) for $(f_i)$ is given by:*

$$A = \inf \left\{ \sum_i \|P_{E^\perp} f_i\|^2; \; E \text{ is a hyperplane in } H \right\},$$

$$B = \sup \left\{ \sum_i \|P_{E^\perp} f_i\|^2; \; E \text{ is a hyperplane in } H \right\}.$$

*Proof.* For any $f \in H$ with $\|f\| = 1$, let $E = [f]^\perp$. Then

$$| < f, f_i > |^2 = \|P_{E^\perp} f_i\|^2, \quad \text{for all } i.$$

It follows immediately that

$$A \geq \inf \left\{ \sum_i \|P_{E^\perp} f_i\|^2; \; E \text{ is a hyperplane in } H \right\},$$

and

$$B \leq \sup \left\{ \sum_i \|P_{E^\perp} f_i\|^2; \; E \text{ is a hyperplane in } H \right\}.$$

Conversely, if $E$ is any hyperplane in $H$, choose $f \in E^\perp$ with $\|f\| = 1$ to get

$$| < f, f_i > |^2 = \|P_{E^\perp} f_i\|^2, \quad \text{for all } i.$$

It follows that

$$A \le \inf \left\{ \sum_i \|P_{E^\perp} f_i\|^2; \ E \text{ is a hyperplane in } H \right\},$$

and

$$B \ge \sup \left\{ \sum_i \|P_{E^\perp} f_i\|^2; \ E \text{ is a hyperplane in } H \right\}.$$

$\square$

An inductive procedure applied to Proposition 4.12 yields the following Corollary.

**Corollary 4.13.** *Let $(f_i)$ be a sequence in a Hilbert space $H$. The following are equivalent:*

*(1) $(f_i)$ is a frame for $H$ with optimal frame bounds $A, B$.*

*(2) For every hyperplane $E \subset H$ we have*

$$A \le \sum_i \|P_{E^\perp} f_i\|^2 \le B,$$

*and A (resp. B) is maximal (resp. minimal) with respect to these two inequalities.*

*(3) For every subspace $E_n \subset H$ of codimension $n$ we have*

$$nA \le \sum_i \|P_{E^\perp} f_i\|^2 \le nB.$$

*and A (resp. B) is maximal (resp. minimal) with respect to these two inequalities.*

The reader should compare this Corollary to Proposition 4.6. Also, give some thought to what this says for normalized tight frames.

A big advantage of treating frames as operators is that theorems about bounded linear operators on a Hilbert space become theorems about frames. For example, a deep question in wavelet theory [138] concerns which families of wavelets are arcwise connected. But the corresponding question for frames has an immediate answer. We just need to recall ([92], page 66) that each pair of partial isometries with the same rank, co-rank and nullity, can be joined by a continuous curve of partial isometries with the same rank, co-rank and nullity (and these two partial isometries cannot be connected without these assumptions). This answers the second part of the following proposition. The first part follows by applying the second part to the polar decomposition of the preframe operator. The **deficiency** of a frame is the dimension of the kernel of the preframe operator (i.e. the codimension of the range of the frame transform).

**Proposition 4.14.** *The family of frames of deficiency $n$ is connected, for all natural numbers $n$. The family of normalized tight frames of deficiency $n$ is connected, for all natural numbers $n$.*

Another example of the power of treating frames as operators is in representing frames as sums of "better objects". This is the underlying idea behind the phase space Wannier functions used in solid state Physics [137, 107, 124] and the Wilson bases [58, 137]. The idea here is to decompose a frame as a sum of two orthonormal bases. In general, Casazza [28] shows that every frame can be written as a sum of three orthonormal bases. This is a consequence of a result from operator theory.

**Theorem 4.15.** *Every bounded operator $T$ on a infinite dimensional complex Hilbert space $H$ can be written in the form $T = a(U_1 + U_2 + U_3)$, where each $U_j$ is a unitary operator and $a$ is a positive real number.*

*Proof.* Fix $0 < \epsilon < 1$ and let

$$S = \frac{1}{2}I + \frac{1-\epsilon}{2}\frac{T}{\|T\|}.$$

Then a calculation shows that $\|I - S\| < 1$, so $S$ is an invertible operator. If we write the polar decomposition of $S$ as $S = VP$, then since $S$ is invertible, $V$ is a unitary operator. Also, we can write (see [92]) $P = \frac{1}{2}(W + W^*)$ where $W$, $W^*$ are unitary. Hence,

$$S = \frac{1}{2}(VW + VW^*),$$

where $VW$, $VW^*$ are unitary. Finally,

$$T = \frac{\|T\|}{1-\epsilon}(VW + VW^* - I),$$

is a sum of three unitary operators. $\qquad\square$

Kalton (see [28]) observed that a bounded, onto, linear operator on $H$ can be written as a linear combination of two unitary operators if and only if it is invertible. In frame language, Theorem 4.15 and Kalton's result become

**Theorem 4.16.** *Every frame (for an infinite-dimensional complex Hilbert space $H$) is (a multiple of) a sum of three orthonormal bases. Moreover, a frame is a linear combination of two orthonormal bases if and only if it is a Riesz basis.*

If we weaken slightly the requirements on the objects we are adding to get a frame, then we can drop down to two. For example [28], every frame is a sum of two normalized tight frames, or an orthonormal basis and a Riesz basis.

As we have seen, if $(f_n)$ is a frame with frame operator $S$, then $(S^{-1}f_n)$ is a frame called the dual frame. Moreover, if we let $g_n = S^{-1}f_n$ then we have by (4.3) that

$$(4.5) \qquad f = \sum_n <f, g_n> f_n, \quad \text{for all } f \in H.$$

In light of (4.5) we define:

**Definition 4.17.** *If $(f_n)$ is a frame for a Hilbert space $H$, a frame $(h_n)$ for $H$ is called an* **alternate dual frame** *(or a* **pseudo-dual***) for $(f_n)$ if*

$$(4.6) \qquad f = \sum_n <f, h_n> f_n, \quad \text{for all } f \in H.$$

We call $(S^{-1}f_n)$ the **canonical dual** of $(f_n)$. If $(f_n)$ is a normalized tight frame, then $S = I$ and so the frame equals its canonical dual. The converse of this clearly holds also. However, in general, a frame may have many dual frames. A simple example would be to let our frame be $\{e_1, e_1, e_2, e_2, \cdots\}$ where $(e_n)$ is an orthonormal basis for $H$ and observe that each of the following is an alternate dual for this frame:

$$\{e_1, 0, e_2, 0, e_3, 0, \cdots\},$$

$$\{0, e_1, 0, e_2, 0, e_3, \cdots\}.$$

The canonical dual for this frame is

$$\{e_1/2, e_1/2, e_2/2, e_2/2 \cdots\}.$$

In Definition 4.17 we assumed that the sequence $(h_n)$ is a frame for $H$. The reason is that there might be sequences satisfying (4.6) which are not frames. For example,

$$\{e_1, \frac{1}{\sqrt{2}}e_2, \frac{1}{\sqrt{2}}e_2, \frac{1}{\sqrt{3}}e_3, \frac{1}{\sqrt{3}}e_3, \frac{1}{\sqrt{3}}e_3, \cdots\}$$

is a normalized tight frame for $H$ and the non-frame sequence

$$\{e_1, \sqrt{2}e_2, 0, \sqrt{3}e_3, 0, 0, \sqrt{4}e_4 0 \cdots\}$$

satisfies (4.6). For the basic properties of alternate dual frames we refer the reader to [93, 103, 108, 126]. Han and Larson [93] have shown that two alternate dual frames are equivalent if and only if they are equal. Also, a frame has a unique alternate dual if and only if it is a Riesz basis. Li [108] has given a characterization of the family of all alternate duals for a given frame. We will return to this topic in the setting of Weyl-Heisenberg frames in the next section.

Later in Section 10 we will return to explore the recent developments in abstract frame theory. For now, we will look at the important specific case of Weyl-Heisenberg frames.

## 5. AN INTRODUCTION TO WEYL-HEISENBERG FRAMES

An excellent "tutorial" on Weyl-Heisenberg frame theory up to 1991 is the paper of Heil and Walnut [96]. A user-friendly introduction to Weyl-Heisenberg frames is the forthcoming book of Gröchenig [88]. The new book edited by Feichtinger and Strohmer [74] is devoted entirely to Gabor (Weyl-Heisenberg) frames. Besides being required reading for anyone who wants to work in this area, [74] also has extensive historical developments concerning all of the notions used here. One should also read the fundamental works of Daubechies [54], and Daubechies, Grossman and Meyer [57]. Also, most of the results here have analogues for $L^2(\mathbb{R}^d)$ for natural numbers $d$. We will not work in this generality, but the basics can be found in a paper of Benedetto [22]. The frames introduced by Gabor [84] are called **Gabor frames** or **Weyl-Heisenberg frames**. The later terminology (introduced in [57]) comes from the representation of the Weyl-Heisenberg group $T \times \mathbb{R} \times \hat{\mathbb{R}}$ acting on $L^2(\mathbb{R})$ by

$$W(x, a, b)f(t) = x \cdot e^{2\pi ib(t-a)}f(t-a).$$

Letting $x = 1$ gives our Weyl-Heisenberg frame.

**Definition 5.1.** *If $a, b \in \mathbb{R}$ and $g \in L^2(\mathbb{R})$ we call $(E_{mb}T_{na}g)_{m,n\in\mathbb{Z}}$ a **Weyl-Heisenberg system** (**WH-system** for short) and denote it by $(g, a, b)$. We call $g$ the **window function**.*

If the WH-system $(g, a, b)$ forms a frame for $L^2(\mathbb{R})$, we call this a **Weyl-Heisenberg frame** (**WH-frame** for short). The numbers $a, b$ are the **frame parameters** with $a$ being the **shift parameter** and $b$ being the **modulation parameter.**

There is a longstanding question concerning WH-frames.

**Problem 5.2.** *Find all $a, b \in \mathbb{R}$ and $g \in L^2(\mathbb{R})$ so that $(g, a, b)$ forms a frame for $L^2(\mathbb{R})$.*

Although WH-frame theory is a very applied area, Problem 5.2 seems to be fundamental to a complete understanding of the field - even if it turns out not to be too useful for specific applications. Later we will see some special cases of this problem which have been solved. There are many variations of this problem which are open. Given a function $g$, what can be said about the set of $a, b \in \mathbb{R}$ for which $(g, a, b)$ have a finite upper frame bound, or forms a frame?

**Problem 5.3.** *Identify those $(g, a, b)$ which have finite upper frame bounds.*

The following problem is interesting because of its relationship to the extended zero divisor conjecture for the case of the Heisenberg group.

**Problem 5.4.** *Given $g \in L^2(\mathbb{R})$, $g \neq 0$, and any finite set $\Lambda \subset \mathbb{R} \times \mathbb{R}$, is the set $(E_b T_a g)_{(a,b) \in \Lambda}$ linearly independent?*

Using $C^*$-algebras, Linnell [110] has shown that Problem 5.4 has a positive answer for $\Lambda \subset a\mathbb{Z} \times b\mathbb{Z}$, for any a,b. Also, Heil, Ramanthan, and Topiwala [95] use the Ergodic Theorem to show that there is a positive answer for any set containing 3 elements as well as a host of other ralated results.

It can be shown by direct calculation that the frame operator $S$ for a WH-frame $(g, a, b)$ commutes with translation by a and modulation by b. Namely, just inner product both sides of the equality in Proposition 5.5 below with an element $h \in L^2(\mathbb{R})$ and simplify.

**Proposition 5.5.** *If $(g, a, b)$ is a WH-frame with frame operator $S$, then for all $h \in L^2(\mathbb{R})$ we have:*

$$S(E_{mb}T_{na}h) = E_{mb}T_{na}Sh, \quad S^{-1}(E_{mb}T_{na}h) = E_{mb}T_{na}S^{-1}h.$$

*In particular, the canonical dual frame of a WH-frame is another WH-frame.*

Another question which could simplify some arguements is

**Problem 5.6.** *Find necessary and sufficient conditions for two WH-frames $(E_{mb}T_{na}g)$ and $(E_{md}T_{nc}h)$ to be equivalent.*

The author asked this question recently at a meeting and in a wonderful display of speed, power and finess, Balan and Landau [10] answered it in a day. They first showed that if two WH-frames $(g, a, b)$ and $(h, c, d)$ are equivalent, then $ab = cd$ and one can deduce from this that the only case that needs to be handled is the case where $b = d$ and $a = c$. Their theorem then states:

**Theorem 5.7.** *Let $(g, a, b)$ and $(h, a, b)$ be two WH-frames. The following are equivalent:*

*(1) $(g, a, b)$ is equivalent to $(h, a, b)$.*
*(2) The linear spans of $(g, \frac{1}{b}, \frac{1}{a})$ and $(h, \frac{1}{b}, \frac{1}{a})$ are equal.*
*(3) If $S^{-1}g$ is the dual frame generator for $(g, a, b)$ then*

$$\|h\|^2 = \sum_m \int h(x)\overline{h(x+ma)} \sum_n \overline{g(x-n/b)}(S^{-1}g)(x-(n/b)-ma).$$

Because this characterization uses the dual frame generator, it is difficult to apply until one answers the following important problem:

**Problem 5.8.** *Given a WH-frame $(g, a, b)$, give an explicit representation of the dual frame generator $S^{-1}g$.*

We point out that Janssen [102] has given a representation for $S^{-1}g$. If $H$ is a matrix with entries $H_{k,\ell;k',\ell'}$ given by

$$< E_{\ell b}T_{ka}g, E_{\ell' b}T_{k'a}g >,$$

Then $(g, a, b)$ is a WH-frame if and only if $H$ represents an invertible operator on $\ell^2(\mathbb{Z}^2)$. Moreover, in this case,

$$S^{-1}(g) = ab \sum_{k,\ell \in \mathbb{Z}} (H^{-1})_{k,\ell;0,0} E_{k/b} T_{\ell/a} g.$$

Finally, Theorem 4.16 states that every WH-frame can be written as a sum of three orthonormal bases (two if $ab = 1$, see below). The problem with Theorem 4.16 is that it uses heavy decomposition results from operator theory, and hence in practice is often not usable. So we ask:

**Problem 5.9.** *For a Weyl-Heisenberg frame $(g, a, b)$, give an explicit representation of this frame as a sum of three orthonormal bases.*

We will be interested in when there are finite upper frame bounds for a WH-system. We call this class of functions the **preframe functions** and denote this class by **PF**. We now have

**Proposition 5.10.** *The following are equivalent:*
*(1) $g \in \mathbf{PF}$.*
*(2) The operator*

$$Sf = \sum_n <f, E_{mb} T_{na} g> E_{mb} T_{na} g,$$

*is a well defined bounded linear operator on $L^2(\mathbb{R})$.*

*Proof.* If $(e_{nm})_{n,m \in \mathbb{Z}}$ is an orthonormal basis for $L^2(\mathbb{R})$, let $T e_{nm} = E_{mb} T_{na} g$. Then $S = TT^*$, and (1) is equivalent to $T$ being a bounded linear operator which, in turn, is equivalent to $TT^*$ being a bounded linear operator. $\square$

For the rest of our work, we will be performing sums of products of certain translates of functions. We will now observe that these always exist and ignore this convergence question hereafter.

**Proposition 5.11.** *If $f, g \in L^2(\mathbb{R})$, and $a, b \in \mathbb{R}$ then for all $k \in \mathbb{Z}$ the series*

$$\sum_{n \in \mathbb{Z}} f(t - na) \overline{g(t - na - k/b)}$$

*converges absolutely a.e. $t \in \mathbb{R}$.*

*Proof.* Since $f, T_{k/b} g \in L^2(\mathbb{R})$, we have that $f \cdot \overline{T_{k/b} g} \in L^1(\mathbb{R})$. Hence,

$$\|f \cdot \overline{T_{k/b} g}\|_{L^1} = \int_{\mathbb{R}} |f(t) \overline{(T_{k/b} g)(t)}| dt = \sum_{n \in \mathbb{Z}} \int_0^a f(t - na) \overline{g(t - na - k/b)} | dt < \infty. =$$

$$\int_0^a \sum_{n \in \mathbb{Z}} |f(t - na) \overline{g(t - na - k/b)}| dt < \infty,$$

where the last equality follows from the Monotone Convergence Theorem. This is all we need. $\qquad\square$

Our next goal is to find some reasonable conditions which guarantee that $(g, a, b)$ forms a WH-frame. Daubechies [54] developed a criteria which led Walnut [134] to use Wiener Amalgam space criteria to produce WH-frames. However, we will first present a much stronger criteria given recently by Casazza and Christensen [36]. To develop this approach, we will make extensive use of the WH-frame identity due to Daubechies [54]. Our proof is due to Heil and Walnut [96].

**Theorem 5.12.** (**WH-Frame Identity.**) *If* $\sum_n |g(t - na)|^2 \leq B$ *a.e. and* $f \in L^2(\mathbb{R})$ *is bounded and compactly supported, then*

$$\sum_{n,m\in\mathbb{Z}} | < f, E_{mb}T_{na}g > |^2 = F_1(f) + F_2(f),$$

*where*

$$F_1(f) = b^{-1} \int_R |f(t)|^2 \sum_n |g(t - na)|^2 dt,$$

*and*

$$F_2(f) = b^{-1} \sum_{k\neq 0} \int_R \overline{f(t)} f(t - k/b) \sum_n g(t - na)\overline{g(t - na - k/b)} dt =$$

$$b^{-1} \sum_{k\geq 1} 2Re \int_R \overline{f(t)} f(t - k/b) \sum_n g(t - na)\overline{g(t - na - k/b)} dt.$$

*Proof.* We are assuming that $f$ is bounded and compactly supported so that all the summations, integrals and interchanges of these below are justified. We define

$$H_n(t) = \sum_k f(t - k/b)\overline{g(t - na - k/b)}.$$

Now, $H_n$ is $1/b$-periodic, $H_n \in L^2[0, 1/b]$ and

$$\int_R f \cdot \overline{E_{mb}T_{na}g(t)} dt = \int_R f(t)\overline{g(t - na)}e^{-2\pi imbt} dt = \int_0^{1/b} H_n(t)e^{-2\pi imbt} dt.$$

Since $(b^{1/2}E_{mb})_{m\in\mathbb{Z}}$ is an orthonormal basis for $L^2[0, 1/b]$, the Plancherel formula yields

$$\sum_m |\int_0^{1/b} H_n(t)e^{-2\pi imbt} dt|^2 = b^{-1} \int_0^{1/b} |H_n(t)|^2 dt.$$

Now we compute

$$\sum_n \sum_m | < f, E_{mb} T_{na} g > |^2 = \sum_n \sum_m | \int_R f(t) \overline{g(t - na)} e^{-2\pi imbt} dt |^2$$

$$= b^{-1} \sum_n \int_0^{1/b} | \sum_k f(t - k/b) \overline{g(t - na - k/b)} |^2 dt$$

$$= b^{-1} \sum_n \int_0^{1/b} \sum_\ell \overline{f(t - \ell/b)} g(t - na - \ell/b) \cdot \sum_k f(t - k/b) \overline{g(t - na - k/b)} dt$$

$$= b^{-1} \sum_n \sum_\ell \int_0^{1/b} \overline{f(t - \ell/b)} g(t - na - \ell/b) \cdot \sum_k f(t - k/b) \overline{g(t - na - k/b)} dt$$

$$= b^{-1} \sum_n \int_R \overline{f(t)} g(t - na) \cdot \sum_k f(t - k/b) \overline{g(t - na - k/b)} dt$$

$$= b^{-1} \sum_k \int_R \overline{f(t)} f(t - k/b) \cdot \sum_n g(t - na) \overline{g(t - na - k/b)} dt$$

$$= b^{-1} \int_R |f(t)|^2 \cdot \sum_n |g(t - na)|^2 dt \ +$$

$$b^{-1} \sum_{k \neq 0} \int_R \overline{f(t)} f(t - k/b) \cdot \sum_n g(t - na) \overline{g(t - na - k/b)} dt.$$

This completes the first part of the WH-Frame Identity. The equality in the last line follows by a simple change of variables. $\qquad \square$

To simplify the notation a little, we introduce the following auxilliary functions:

$$(5.1) \qquad G_k(t) = \sum_{n \in \mathbb{Z}} g(t - na) \overline{g(t - na - k/b)}, \quad \text{for all } k \in \mathbb{Z}.$$

Note that the $G_k$ are periodic functions on $\mathbb{R}$ of period a. During our study of WH-frames, the reader will become aware that the $G_k$ above contain most of the important information about a Weyl-Heisenberg frame. The trick is to see how it applies to a given situation.

There are several important consequences of the WH-Frame Identity which we now examine. The first consequence of the WH-frame identity comes from [38].

**Corollary 5.13.** *Let $a, b \in \mathbb{R}$ with $ab \leq 1$ and $g \in L^2(\mathbb{R})$ and assume that*

$$\sum_k |G_k(t)|^2 \leq B, \quad a.e.$$

*Then for all bounded, compactly supported functions $f \in L^2(\mathbb{R})$ the series*

$$Lf = b^{-1} \sum_k (T_{k/b}f)G_k,$$

*converges unconditionally in norm in $L^2(\mathbb{R})$. Moreover,*

$$< Lf, f > = \sum_{m,n \in \mathbb{Z}} | < f, E_{mb}T_{na}g > |^2.$$

*Finally, if $g \in \mathbf{PF}$, so that the series*

$$Sf = \sum_{m,n \in \mathbb{Z}} < f, E_{mb}T_{na}g > E_{mb}T_{na}g$$

*also converges unconditionally in $L^2(\mathbb{R})$, we have that $Lf = Sf$.*

*Proof.* First we check that the series for $Lf$ converges unconditionally for all bounded compactly supported $f \in L^2(\mathbb{R})$. But these functions are finite sums of bounded functions supported on intervals of the form $I_n = [na, (n+1)a]$. So we assume $f$ is supported on $I_n$ with uniform upper bound $D$. Now, since $a \leq 1/b$, we have that the functions $((T_{k/b}f)G_k)_{k \in \mathbb{Z}}$ are disjointly supported. Since the $G_k$ are periodic of period $a$, a simple calculation yields

$$\| \sum_{k \in M} (T_{k/b}f)G_k \|_{L^2(\mathbb{R})} = \int_0^a |f(t)|^2 \sum_{k \in M} |G_k(t)|^2 dt \leq D^2 \int_0^a \sum_{k \in M} |G_k(t)|^2 dt.$$

Since $\sum_k |G_k(t)|^2 \leq B$ a.e., an application of the Monotone Convergence Theorem yields that our series (even every subseries) converges. So the series for $Lf$ converges unconditionally in $L^2(\mathbb{R})$.

For the moreover part, we check

$$< Lf, f > = < b^{-1} \sum_{k \in \mathbb{Z}} (T_{k/b}f)G_k, f > = b^{-1} \sum_{k \in \mathbb{Z}} < (T_{k/b}f)G_k, f >$$

$$= b^{-1} \sum_{k \in \mathbb{Z}} \int_R \overline{f(t)} f(t - k/b)G_k(t)dt = \sum_{m,n \in Z} | < f, E_{mb}T_{na}g > |^2,$$

where the last equality follows from the WH-Frame Identity.

To see that $Sf = Lf$, we redo the proof of the WH-Frame Identity to see that for all $h \in L^2(\mathbb{R})$ we have

$$< Sf, h > = \sum_{m,n \in \mathbb{Z}} < f, E_{mb}T_{na}g > < E_{mb}T_{na}g, h > =$$

$$< b^{-1} \sum_k (T_{k/b} f) \cdot G_k, h > .$$

$\square$

Our next consequence of the WH-Frame Identity is a necessary condition for $(g, a, b)$ to generate a WH-frame.

**Corollary 5.14.** *If $(g, a, b)$ generates a WH-frame with frame bounds $A, B$ then*

$$Ab \leq G_0(t) = \sum_n |g(t - na)|^2 \leq Bb, \quad a.e.$$

*In particular, $g$ must be bounded.*

*Proof.* For any bounded function $f$ supported on an interval $I$ of length $\leq 1/b$ we have that $F_2(f) = 0$ in the WH-Frame Identity and so

$$A\|f\|^2 \leq \sum_{n,m \in \mathbb{Z}} | < f, E_{mb} T_{na} g > |^2 = b^{-1} \int_R |f(t)|^2 G_0(t) dt \leq B\|f\|^2.$$

The result follows easily from here. $\square$

It can be shown that the condition $A \leq G_0(x) \leq B$ a.e. is equivalent to $(E_{m/a} g)$ is a Riesz basic sequence (i.e. a Riesz basis for its span) (see Theorem 9.2 below).

Since the Fourier transform is a unitary operator on $L^2(\mathbb{R})$ which takes $E_{mb} T_{na} g$ to $T_{mb} E_{-na} \hat{g}$, it follows that $(g, a, b)$ generates a WH-frame if and only if $(\hat{g}, b, a)$ generates a WH-frame. Hence, by Corollary 5.14 both $g$ and $\hat{g}$ must be bounded functions. Although the conditions in Corollary 5.14 are considered basic for WH-frames, Casazza and Christensen [36] studied WH-systems which are frames for subspaces of $L^2(\mathbb{R})$ and showed that in this case, the lower inequality is no longer necessary for the existence of a frame sequence. Gabardo and Han [82] have generalized these results to Weyl-Heisenberg unitary systems.

**Corollary 5.15.** *If $g \in L^2(\mathbb{R})$, $a, b \in \mathbb{R}$ and supp $g \subset I \subset R$, where $I$ is an interval of length $\leq 1/b$, then $(g, a, b)$ forms a WH-frame if and only if there are constants $A, B > 0$ so that*

$$A \leq G_0(t) \leq B, \quad a.e.$$

*Proof.* Since $g$ is supported on an interval of length $1/b$, we have that $G_k(t) = 0$ a.e. for all $k \neq 0$. Now, by the WH-Frame Identity we have that $F_2(f) = 0$ for all bounded, compactly supported $f \in L^2(\mathbb{R})$ and so

$$\sum_{m,n \in \mathbb{Z}} | < f, E_{mb} T_{na} g > |^2 = b^{-1} \int_R |f(t)|^2 G_0(t) dt.$$

It follows that we have the frame inequalities (3.1) for all bounded, compactly supported $f \in L^2(\mathbb{R})$. It is an instructive problem to show that a squence $(f_n)$ satisfying the frame inequalities for a dense set of vectors in $H$ must satisfy these inequalities on the whole space $H$. □

We give one more consequence of the WH-frame Identity.

**Corollary 5.16.** *Let $a, b \in \mathbb{R}$ and assume that $g$ is supported on an interval $I$ with $|I| \leq \frac{1}{a}$. If*

$$\sum_{m \in \mathbb{Z}} |g(t - nb)|^2 = a, \quad a.e.$$

*then $(E_{mb}T_{na}\hat{g})$ is a normalized tight WH-frame consisting of non-compactly supported functions.*

*Proof.* From the WH-frame identity, one can easily check that $(E_{ma}T_{nb}g)$ is a normalized tight WH-frame for $L^2(\mathbb{R})$. Taking the Fourier transform of this yields the result. □

We have just seen that there are some restrictions on which functions $g \in L^2(\mathbb{R})$ can form WH-frames. Another restriction is given by the important Balian-Low Theorem (proved independently by Balian [9] and Low [111]). The original proofs contained a technical gap which was later filled in by Coifman and Semmes (see [54] for this corrected proof). There is also a proof of this theorem by Battle [12] for the orthonormal basis case based on the Heisenberg Uncertainty Principle. There are other related proofs due to Daubechies and Janssen [59] and by Benedetto, Heil and Walnut [17].

**Theorem 5.17** (Balian-Low). *If $g \in L^2(\mathbb{R})$, $ab = 1$ and $(g, a, b)$ generates a WH-frame, then either $tg(t) \notin L^2(\mathbb{R})$ or $g' \notin L^2(\mathbb{R})$.*

It follows that if $ab = 1$, and $(g, a, b)$ is a WH-frame for $L^2(\mathbb{R})$, then either $g$ does not decay very rapidly, or $g$ is not smooth. In particular, the Gaussian functions used by Gabor in his original work cannot yield frames when $ab = 1$, despite the fact that this WH-system is complete in $L^2(\mathbb{R})$.

It is instructive to compare the Balian-Low Theorem to the Classical Uncertainty Principle Inequaltiy. The elegent proof below is due to Wiener and can most easily be found in Benedetto, Heil and Walnut [17].

**Theorem 5.18. (Classical Uncertainty Principle Inequality.)** *Let $(t_0, \gamma_0) \in \mathbb{R} \times \mathbb{R}$. Then for every $f \in L^2(\mathbb{R})$ we have*

$$(5.2) \qquad \|f\|_2^2 \leq 4\pi \|(t - t_0)f(t)\|_2 \|(\gamma - \gamma_0)\hat{f}(\gamma)\|_2.$$

*Moreover, there is equality in (5.2) if and only if $f(t) = CE_{\gamma_0}e^{-s(t-t_0)^2}$, for some complex number $C$ and $s > 0$.*

*Proof.* The mapping $f(t) \mapsto E_{-\gamma_0} T_{-t_0} f(t)$ shows that we only need to verify (5.2) for $(t_0, \gamma_0) = (0, 0)$. Now we just compute for this case,

$$\|f\|_2^4 = \left( \int_R r|f(t)^2|'dt \right)^2 \leq \left( \int_R |t||f(t)^2|'dt \right)^2$$

$$\leq 4 \left( \int_R |t\overline{f(t)}f'(t)|dt \right)^2 \leq 4\|tf(t)\|_2^2\|f'(t)\|_2^2 = 16\pi^2\|tf(t)\|_2^2\|\gamma\hat{f}(\gamma)\|_2^2.$$

The above works for "nice" $f$. We then extend to all $f$ by standard duality arguements. We leave it for the reader to check the "moreover" part of the Theorem. □

For some signals, inequality (5.2) does not provide useful information. For example, if $f \in L^2(\mathbb{R})$ behaves like $|t|^a$ as $t \to \infty$, where $a \in [-3/2, -1/2)$, then the right hand side of the inequality is infinite. Also, by the Balian-Low Theorem, if $ab = 1$ and $(g, a, b)$ gives a WH-frame, then the right side of (5.2) is infinite if $f$ is replaced by $g$. That is, the Balian-Low Theorem maximizes the Classical Uncertainty Principle Inequality. To relate this to our original piano concert example, inequality (5.2) asserts that if a sound $f$ is emitted at time $t_0$ and lasts a very short time (so $(t - t_0)$ is very small) then the frequency range for $f$ is quite broad. That is, $f$ cannot be very close to a pure tone of frequency $\gamma_0$ (for if it were then both $\|(t-t_0)f(t)\|_2^2$ and $\|(\gamma - \gamma_0)\hat{f}(\gamma)\|_2^2$ would be small in contrast to the "loudness" $\|f\|_2^2$).

The above results put restrictions on $g$ in order that $(g, a, b)$ yield a WH-frame. There are also some restrictions on $a, b$.

**Proposition 5.19.** *Let $g \in L^2(\mathbb{R})$ and $a, b \in \mathbb{R}$.*
  *(1) If $(E_{mb}T_{na}g)$ is complete, then $ab \leq 1$.*
  *(2) If $(g, a, b)$ is a WH-frame and*
    *(i) $ab < 1$ then $(g, a, b)$ is overcomplete.*
    *(ii) $ab = 1$ then $(g, a, b)$ is a Riesz basis.*

Part (1) of Proposition 5.19 has a complicated history (see [54] for a discussion) which derives from the work of Rieffel [125]. Today, there is a simpler proof using Beurling density due to Ramanathan and Steger [123]. Moreover, the results of Remanathan and Steger [123] combined with an important example of Benedetto, Heil and Walnut [17] shows that the form of the lattice in the Rieffel result [125] is quite important to the conclusion. There are many derivations available for (2) [38, 54, 55, 93, 101, 102]. A proof of (2) will be given later in Remark 9.4. For now we ask,

**Problem 5.20.** *Given $g \in L^2(\mathbb{R})$ and $ab > 1$, write down explicitly a function $f \in L^2(\mathbb{R})$ for which $f \perp span\ E_{mb}T_{na}g$. In fact, we would conjecture that if*

$n < ab \leq n + 1$, then there are functions $(f_i)_{i=1}^n$ in $L^2(\mathbb{R})$ so that

$$H = (span\ E_{mb}T_{na}g) \oplus \sum_{i=1}^n \oplus (span\ E_{mb}T_{na}f_i).$$

where the sums above are orthogonal sums.

For $ab = 2$, Daubechies [54] gives the explicit representation for the function $f$ requested in Problem 5.20. Just before this article went to print, Balan and Landau [10] and independently Gabardo and Han [83] answered it positively.

It can be checked by direct calculation that a necessary condition for $(g, a, b)$ to form a WH-frame is that $\sum_k |G_k(t)|^2 \leq B$ a.e. This condition is not sufficient however (see Example 5.35). We are now ready to pass to the CC-condition [36] which is sufficient for $(g, a, b)$ to have a finite upper frame bound.

**Theorem 5.21** (CC-Condition). *If $g \in L^2(\mathbb{R})$, $a, b \in \mathbb{R}$ and*

(CC) $$\sum_{k \in \mathbb{Z}} |\sum_{n \in \mathbb{Z}} g(t - na)\overline{g(t - na - k/b)}| = \sum_{k \in \mathbb{Z}} |G_k(t)| \leq B, \quad a.e.,$$

*then $g \in \mathbf{PF}$. Moreover, if we also have*

(5.3) $$\sum_{k \neq 0} |G_k(t)| \leq (1 - \epsilon)G_0(t) \quad a.e.,$$

*for some $0 < \epsilon < 1$, then $(g, a, b)$ is a WH-frame.*

*Proof.* We first observe that

$$\sum_{k \neq 0} |T_{-k/b}G_k(t)| = \sum_{k \neq 0} |T_{-k/b}\sum_{n \in \mathbb{Z}} T_{na}g(t)\overline{T_{na+k/b}g(t)}|$$

$$= \sum_{k \neq 0} |\sum_{n \in \mathbb{Z}} T_{na-k/b}g(t)\overline{T_{na}g(t)}| = \sum_{k \neq 0} |\sum_{n \in \mathbb{Z}} T_{na+k/b}g(t)\overline{T_{na}g(t)}|$$

$$= \sum_{k \neq 0} |\sum_{n \in \mathbb{Z}} \overline{T_{na+k/b}g(t)}T_{na}g(t)| = \sum_{k \neq 0} |G_k(t)|.$$

Now,

$$|\sum_{k \neq 0} \int \overline{f(t)}f(t - k/b)\sum_{n \in \mathbb{Z}} g(x - na)\overline{g(x - na - k/b)}dt|$$

$$\leq \sum_{k \neq 0} \int |f(t)| \cdot |T_{k/b}f(t)| \cdot |G_k(t)|dt$$

$$\leq \sum_{k \neq 0} \left(\int |f(t)|^2|G_k(t)|dt\right)^{1/2} \cdot \left(\int |T_{k/b}f(t)|^2|G_k(t)|dt\right)^{1/2}$$

$$\leq \left( \sum_{k \neq 0} \int |f(t)|^2 |G_k(t)| dt \right)^{1/2} \cdot \left( \sum_{k \neq 0} \int |T_{k/b} f(t)|^2 |G_k(t)| dt \right)^{1/2}$$

$$= \left( \int |f(t)|^2 \sum_{k \neq 0} |G_k(t)| dt \right)^{1/2} \cdot \left( \int |f(t)|^2 \sum_{k \neq 0} |T_{-k/b} G_k(t)| dt \right)^{1/2}$$

$$= \int |f(t)|^2 \sum_{k \neq 0} |G_k(t)| dt.$$

Combining this with a simple application of the triangle inequality to the WH-Frame Identity yields the result. $\qquad \square$

The first condition in Theorem 5.21 is called the **CC-condition** and plays an important role in Weyl-Heisenberg frame theory [39] which we will explore further in Section 9. It is known (see Example 5.35 below) that the CC-condition is not necessary for having a finite upper frame bound or a WH-frame. A natural open question concerns the CC-condition and the dual frame.

**Problem 5.22.** *If $(g, a, b)$ satisfies the CC-condition, does the canonical dual $(S^{-1}g, a, b)$ also satisfy the CC-condition?*

There is also a **uniform CC-condition** which is slightly stronger than the CC-condition and does pass to the dual frame - at least in the case where $ab$ is rational (see [39], Theorem 4.14). We say that $(g, a, b)$ satisfies the **uniform CC-condition** if for every $\epsilon > 0$ there is a natural number $K > 0$ so that for a.e. t we have

$$\text{(UCC)} \qquad\qquad \sum_{|k| \geq K} |G_k(t)| < \epsilon.$$

There are WH-frames $(g, a, b)$ which satisfy the CC-condition but fail the uniform CC-condition (see example 5.34 below). Even the uniform CC-condition is a fairly weak assumption since the condition of Tolimieri and Orr (see [103], Subsection 1.4.3) is strong enough to imply it (see [39], Proposition 4.12).

Until now, our WH-frames $(g, a, b)$ have been heavily dependent on the values of a and b in the sense that a small change in either one of these values may cause our family to cease to be a frame. Now we will consider some results of Daubechies [54] and Walnut [134] which guarantee that certain functions will form WH-frames for a fixed value of a and all small values of b. To do this we define the **Wiener amalgam space** $W(L^\infty, \ell^1)$ to be the set of all measurable functions $g$ on $\mathbb{R}$ for which there is some $a > 0$ such that

$$(5.4) \qquad \|g\|_{W,a} = \sum_{n \in \mathbb{Z}} \|g \cdot \chi_{[an, a(n+1))}\|_\infty = \sum_{n \in \mathbb{Z}} \|T_{na} g \cdot \chi_{[0,a)}\|_\infty < \infty.$$

It is easily checked that $W(L^\infty, \ell_1)$ is a Banach space with the norm $\|\cdot\|_{W,a}$. Straightforward calculations (see e.g. [96]) yield the following proposition.

**Proposition 5.23.** *For a function $g \in W(L^\infty, \ell_1)$ we have*
  *(1) If $\|g\|_{W,a}$ is finite for one value of $a > 0$, then it is finite for all $a > 0$.*
  *(2) If $m$ is a natural number and $0 < b \le ma$, then $\|g\|_{W,a} \le 2m\|g\|_{W,b}$.*

The following result is a generalization by Walnut [134] of a result of Daubechies [54]. The proof of Daubechies uses the Poisson summation formula. Our proof is due to Heil and Walnut [96].

**Theorem 5.24.** *Let $g \in L^2(\mathbb{R})$ and $a > 0$ be such that:*
  *(1) There exist constants $A, B$ such that $0 < A \le G_0(t) \le B$ a.e. $t \in \mathbb{R}$,*
  *(2) We have $g \in W(L^\infty, \ell_1)$.*
  *Then there exists a $b_0 > 0$ so that $(g, a, b)$ is a WH-frame for all $0 < b \le b_0$.*

To prove Theorem 5.24, one just checks that the CC-condition and the equation following it holds for all small values of $b$ (see e.g., [96], Theorem 4.1.5). Some important examples related to Theorem 5.24 are due to Feichtinger and Janssen [78]. They show that WH-frame bound conditions depend heavily on the lattice parameters. They show in particular that the family of parameters $a, b$ for which $(g, a, b)$ forms a WH-frame need not be intervals, and can vary greatly for rational vrs irrational parameters.

We end this section by introducing the Zak transform (also known as the Weil-Brezin map). This is one of the most powerful tools in WH-frame theory for producing examples, especially in the case $ab = 1$. This operator was introduced by Zak [140] in the 1960's (calling it the **kq-representation**) when he was working in solid state physics. However, the transform was around prior to this (see the introduction of [74] for a detailed account of the history). Since then, it has been extensively used in WH-frame theory, especially by Janssen (alone and with various co-authors) [23, 99, 100, 38, 39]. A good introduction to the Zak transform is [100]. For extensive examples using the Zak transform see [39]. We start with the operator theoretic definition for the case $a = b = 1$. In this case, $(E_m T_n \chi_{[0,1]})$ forms an orthonormal basis for $L^2(\mathbb{R})$.

**Definition 5.25.** *The **Zak Transform** is the map $Z : L^2(\mathbb{R}) \to L^2(\mathbb{Q})$ ($\mathbb{Q} = [0,1) \times [0,1)$) given by*
$$Z(E_m T_n \chi_{[0,1]}) = E_{m,n}.$$

Normally, the following definition is used for the Zak Transform. As we will see afterwards, these two definitions are equivalent.

**Definition 5.26.** *For $\lambda > 0$, the **Zak transform** of a function $f \in L^2(\mathbb{R})$ is*

$$(5.5) \qquad (Z_\lambda f)(t, \nu) = \lambda^{1/2} \sum_{k \in \mathbb{Z}} f(\lambda(t - k)) e^{2\pi i k \nu}, \quad a.e.\ t, \nu \in \mathbb{R},$$

*where the right-hand side has to be interpreted in the $L^2_{loc}(\mathbb{R})$ sense. If no confusion will arise, we will write $Z$ for $Z_\lambda$.*

We check that the series in the Zak transform converges for the case $\lambda = 1$ (the general case follows with only notational changes). We choose $f \in L^2(\mathbb{R})$ and $k \in \mathbb{Z}$ and set

$$M_k(t, \nu) = f(t - k)e^{2\pi i k \nu}.$$

Now we compute

$$\|M_k\|_2^2 = \int_0^1 \int_0^1 |f(t-k)e^{2\pi i k \nu}|^2 d\nu dt = \int_0^1 |f(t-k)|^2 dt < \infty.$$

So $M_k \in L^2(Q)$. Moreover, these functions are orthogonal. That is, for all $j \neq k$

$$< M_k, M_j > = \int_0^1 f(t-j)\overline{f(t-k)} \left( \int_0^1 e^{2\pi i(j-k)\nu} d\nu \right) dt = 0.$$

It follows that $\| \sum_k M_k \|_2^2 = \sum \|M_k\|_2^2 = \|f\|_2^2$. Therefore, $\sum_k M_k$ is well defined, linear, convergent, and norm-preserving.

From the definition, we have the **quasi-periodicity relations**

$$(5.6) \qquad Zf(t+1, \nu) = e^{-2\pi i \nu} Zf(t, \nu) \quad \text{and} \quad Zf(t, \nu+1) = Zf(t, \nu).$$

We see then that the Zak transform is completely determined by its values in the unit square $Q = [0, 1) \times [0, 1)$. So we can define a new Hilbert space

$$L^2(Q) = \{F : Q \to C : \|F\|_2 = \left( \int_0^1 \int_0^1 |F(t, \nu)|^2 d\nu dt \right)^{1/2} < \infty\}.$$

The inner product here is

$$< F, G > = \int_0^1 \int_0^1 F(t, \nu)\overline{G(t, \nu)} d\nu dt.$$

Note that the two dimensional exponentials $(E_{mn})_{m,n \in \mathbb{Z}}$ form an orthonormal basis for $L^2(Q)$. Now we are ready to relate the Zak transform to our Hilbert spaces.

**Theorem 5.27.** *The Zak transform is a unitary map of $L^2(\mathbb{R})$ onto $L^2(Q)$.*

*Proof.* Again we will check the case $\lambda = 1$. For $m, n \in \mathbb{Z}$, let $\phi_{mn} = T_n E_m \chi_{[0,1)}(t)$. A direct calculation shows that $(\phi_{mn})$ is an orthonormal basis for $L^2(\mathbb{R})$. We will show that $Z$ maps $\phi_{mn}$ to $E_{mn}$ to complete the proof. Now,

$$Z\phi_{mn}(t, \nu) = \sum_{k \in \mathbb{Z}} e^{2\pi i m(t+k-n)} \chi_{[n,n+1)}(t+k)e^{2\pi i k \nu}.$$

Since $0 \leq t \leq 1$, the only nonzero term in this series is $k = n$, so

$$Z\phi_{mn}(t, \nu) = e^{2\pi i m t} e^{2\pi i n \nu} = E_{mn}(t, \nu).$$

$\square$

We can also recover a function from its Zak transform by the formula

$$(5.7) \qquad\qquad f(t) = \int_0^1 (Zf)(t, \nu) d\nu.$$

A close look at the proof of Theorem 5.27 yields for $ab = 1$ that

$$(5.8) \qquad\qquad Z(E_{mb} T_{na} g) = E_{mn} Z g.$$

Equation (5.8) places some strong restrictions on the form of $Zg$ when $(g, a, b)$ yields a frame. We combine some of these results in the following theorem [54, 94, 99].

**Theorem 5.28.** *Let $ab = 1$ and $g \in L^2(\mathbb{R})$.*
*(1) $(E_{mb} T_{na} g)$ is complete in $L^2(\mathbb{R})$ if and only if $Zg \neq 0$ a.e.*
*(2) $(g, a, b)$ generates an orthonormal basis for $L^2(\mathbb{R})$ if and only if $|Zg| = 1$ a.e.*
*(3) The following are equivalent:*
*(i) $0 < A \leq |Zg|^2 \leq B$ a.e.*
*(ii) $(g, a, b)$ generates a frame for $L^2(\mathbb{R})$ with frame bounds $A, B$.*
*(4) $(g, a, b)$ is minimal if and only if $(1/Zg) \in L^2(\mathbb{Q})$.*

*Proof.* (1). This is immediate from the fact that the Zak transform is a unitary map, $(E_{mn})$ is an orthonormal basis for $L^2(Q)$, and (5.8).

(2). We leave this to the reader.

(3). Since $Z$ is a unitary map, (5.8) yields that $(g, a, b)$ is a frame for $L^2(\mathbb{R})$ with frame bounds $A, B$ if and only if $(E_{mn} Z g)$ is a frame for $L^2(Q)$ with frame bounds $A, B$. But, for any $F \in L^2(Q)$,

$$\sum_{m,n} |< F, E_{mn} Z g >|^2 = \sum_{m,n} |< F \cdot \overline{Zg}, E_{mn} >|^2 = \| F \cdot Zg \|^2.$$

The equivalence of $(i)$ and $(ii)$ follows easily from here. Finally, (4) is an observation we leave to the reader. $\square$

We mention yet another important result which has a long history (see [99]).

**Theorem 5.29.** *Let $f \in L^2(\mathbb{R})$ and assume that $Zf$ is continuous on $\mathbb{R} \times \mathbb{R}$. Then $Zf$ has a zero.*

Combining Theorem 5.28 with Theorem 5.29 we see that no function whose Zak transform is continuous can generate a WH-frame for $ab = 1$. However, if we are willing to go to $ab < 1$ then we can get good functions, even Gaussians, to be window functions.

**Remark 5.30.** *It can be checked that the Zak transform of the Gaussian function $g(t) = \pi^{-1/4}e^{-t^2/2}$ is continuous and has a single zero in $Q$. It follows from Theorem 5.28 and the Balian-Low Theorem that (when $ab = 1$) $(E_{mb}T_{na}g)$ is complete in $L^2(\mathbb{R})$ but does not form a frame. However, it does form a frame for all values of $ab < 1$. This was shown independently by Lyubarskii [112] and Seip and Wallsten [131]. Recently, and Lyubarskii and Seip [114] proved that removal of just one element from this WH-system leaves an independent set. Also, Casazza and Lammers [42] have written down explicitly for $ab = 1$ those functions $g \in L^2(\mathbb{R})$ for which $(g, a, b)$ is complete.*

This raises another question for WH-systems.

**Problem 5.31.** *Find all those $g, a, b$ so that $(g, a, b)$ is complete in $L^2(\mathbb{R})$.*

Problem 5.31 has been answered for the case $ab = 1$ by Casazza and Lammers [42]. But the general case is completely open. We now state another important result for applications of the Zak transform.

**Proposition 5.32.** *Suppose that $g \in \mathbf{PF}$. For each $k \in \mathbb{Z}$ we have*
*(1) $|(Zg)(t, \nu)|^2 \leq B$, a.e. $t, \nu \in [0, 1)$.*
*(2) $|(Zg)(t, \nu)|^2$ has for a.e. $t$ the Fourier series expansion*

$$|(Zg)(t, \nu)|^2 = \sum_{k=-\infty}^{\infty} G_k(t)e^{-2\pi i k\nu}, \quad a.e. \ \nu.$$

*(3) For a.e. $t$ and all $k \in \mathbb{Z}$,*

$$G_k(t) = \int_0^1 |(Zg)(t, \nu)|^2 e^{2\pi i k\nu}d\nu.$$

*Proof.* (1) follows from the definition of the Zak transform. (2) is an application of Parseval's identity and Carleson's theorem, and (3) follows immediately from (2).  □

We end our introduction to WH-frames by considering the basic use of the Zak transform to produce examples in this area. This result can be found in [39].

**Proposition 5.33.** *Let $a = b = 1$ and $g \in \mathbf{PF}$. For $M \subset Z$ with $|M| < \infty$ and all $f \in L^2(\mathbb{R})$ let $S_M : L^2(\mathbb{R}) \to L^2(\mathbb{R})$ be given by*

$$S_M f = \sum_{k \in M} f(\cdot - k)G_k = \sum_{k \in M}(T_k f) \cdot G_k.$$

*Then*

$$\|S_M\| = ess\ sup_{t,\nu} |\sum_{k \in M} G_k(t) e^{-2\pi i k \nu}|.$$

*Moreover, for all $\nu_0 \in [0, 1)$ we have*

$$ess\ sup_t |\sum_{k \in M} G_k(t) e^{-2\pi i k \nu_0}| \leq \|S_M\|.$$

*Proof.* For any $f \in L^2(\mathbb{R})$, taking the Zak transform, we have

$$Z\left[\sum_{k \in M} f(\cdot - k) G_k\right](t, \nu) = (Zf)(t, \nu) \cdot \sum_{k \in M} G_k(t) e^{-2\pi i k \nu}, \quad \text{a.e.} \quad t, \nu \in \mathbb{R}.$$

Since $Z$ is unitary, the operator norm of $S_M$ is the same as the operator norm of the multiplication operator:

$$Zf \in L^2(Q) \rightarrow (Zf)(t, \nu) \cdot \sum_{k \in M} G_k(t) e^{-2\pi i k \nu}.$$

But these operator norms are precisely

$$\|S_M\| = ess\ sup_{t,\nu} |\sum_{k \in M} G_k(t) e^{-2\pi i k \nu}|.$$

We leave the moreover part to the reader. $\qquad\square$

Proposition 5.33 is the main tool for constructing counterexamples using the Zak transform. We will give two such examples here to give the flavor of this method. These examples come from [39].

**Example 5.34.** *There is a WH-frame $(g, 1, 1)$ satisfying the CC-condition but failing the uniform CC-condition.*

*Proof.* Choose real numbers $a_m > 0$ with $M = \sum_m a_m < \infty$. Define $F(t, \nu)$ for $t, \nu \in [0, 1)$ by

$$F(t, \nu) = M + \sum_{m=\ell}^{\infty} a_{m-\ell} cos\ 2\pi m \nu, \quad t \in [1 - \frac{1}{\ell}, 1 - \frac{1}{\ell+1}), \ \nu \in [0, 1).$$

Let $g \in L^2(\mathbb{R})$ be the unique element satisfying

$$(Zg)(t, \nu) = \sqrt{F(t, \nu)}, \ t, \nu \in [0, 1).$$

A direct calculation for this $g$ yields that $G_k = 0$ for $|k| = 0, 1, \cdots m - 1$ while $G_k(t) = \frac{1}{2} a_{|k|-\ell}$ for $|k| = \ell, \ell+1, \cdots$, when $t \in [1 - 1/\ell, 1 - 1/(\ell + 1)$. A moment's reflection should convince the reader that this implies that $(g, 1, 1)$ satisfies the CC-condition while failing the uniform CC-condition. $\qquad\square$

**Example 5.35.** *There is a function $g \in L^2(\mathbb{R})$ so that $(g, 1, 1)$ yields a WH-frame for $L^2(\mathbb{R})$, but $g$ fails to satisfy the CC-condition.*

*Proof.* We see that what we need is a function $Zg$ so that $Zg$ is essentially bounded but $(Zg)(t, \cdot)$ is not continuous a.e. t. So let

$$(Zg)(t, \nu) = F(t, \nu) = \begin{cases} 1: & 0 \leq t \leq 1, \quad 0 \leq \nu \leq \frac{1}{2} \\ \frac{1}{2}: & 0 \leq t < 1, \quad \frac{1}{2} < \nu < 1 \end{cases}$$

By Theorem 5.28 we see that $(g, 1, 1)$ yields a WH-frame for $L^2(\mathbb{R})$. Using Proposition 5.32 we can explicitly compute:

$$G_k(t) = \frac{3((-1)^k - 1)}{8\pi ik}, \quad k \neq 0.$$

It follows that

$$\sum_k |G_k(t)| = \infty.$$

So our example is complete. $\qquad\qquad\qquad\qquad\qquad\qquad\qquad\qquad\square$

## 6. An Introduction to Perturbation Theory

Perturbation theory involves answering the question: If $(f_i)$ is a frame for a Hilbert space $H$ and $(g_i)$ is a sequence in $H$ which is "close" to $(f_i)$, must $(g_i)$ be a frame for $H$ (or just for its span)? And if so, does it have the same excess? What can be said about $(g_i)$ if $(f_i)$ is a Riesz basis? The main issues here are (1) How do we measure "closeness"? and (2) Once we have a measure of "closeness" how close does the sequence need to be to a frame to guarantee that it is a frame? For the first question, there are two natural ways to measure closeness in a Hilbert space: Close *in norm* and close *in inner product*. We will look at both these possibilities here and also look at natural conditions for question (2). Again, the important book of Young [139] contains many results and references for perturbation theory in the context of frames - especially Fourier frames which we consider in Section 8.

Perhaps the first perturbation theorem was due to C. Neumann [119]. This states that an operator $T$ on a (Banach) Hilbert space $H$ satisfying $\|I - T\| < 1$, must be an invertible operator. As a consequence we get a perturbation theorem:

**Theorem 6.1.** *If $(f_i)$ and $(g_i)$ are sequences in a Hilbert space $H$ and $0 < \lambda < 1$ satisfies for all sequences of scalars $(a_i)$,*

$$\|\sum_i a_i(f_i - g_i)\| \leq \lambda \|\sum_i a_i f_i\|,$$

*then $Tf_i = g_i$ is a well-defined invertible operator on $H$.*

*Proof.* Clearly, $\|I - T\| \leq \lambda < 1$.                                        $\square$

Theorem 6.1 is generally referred to as the Paley-Wiener Perturbation Theorem and is proved in [121]. Apparently no one noticed that it was really due to Carl Neumann.

A much more powerful perturbation theorem is due to Hilding [97] and is a generalization of the perturbation theorem of C. Neumann (See also [43] for a generalization) and is also really a Banach space result.

**Theorem 6.2.** *Let* $(f_i)$ *and* $(g_i)$ *be sequences in a (Banach) Hilbert space* $H$. *Assume there are constants* $\lambda_1, \lambda_2$ *satisfying for all sequences of scalars* $(a_i)$

$$\|\sum_i a_i(f_i - g_i)\| \leq \lambda_1 \|\sum_i a_i f_i\| + \lambda_2 \|\sum_i a_i g_i\|,$$

*then* $Tf_i = g_i$ *is a well-defined invertible operator on* $H$.

The previous two results say that if $(f_i)$ is a frame for $H$ and if $(g_i)$ is a sequence in $H$ which is close to $(f_i)$ in the sense of the inequalities in the theorems, then $(g_i)$ is a frame for $H$ *which is equivalent to the frame* $(f_i)$.

It would be advantageous to be able to conclude that the $(g_i)$ forms a frame - *even if it is not equivalent to the frame* $(f_i)$. Such results existed in exponential frame theory for quite some time (see Section 8). We next give a perturbation theorem which is an unpublished observation of Casazza. The advantage of this result is that it is the first "necessary and sufficient" perturbation theorem, and hence is best possible in all situations. Also, it allows one to conclude that our second family forms a frame without it being equivalent to the original sequence. Finally, it is easy to reduce most abstract perturbation theorems to this theorem quite quickly.

**Theorem 6.3.** *Let* $(f_i)$ *be a frame for a Hilbert space* $H$ *and let* $(g_i)$ *be a sequence of elements of* $H$. *The following are equivalent:*

*(1)* $(g_i)$ *is a frame for* $H$.
*(2) There is a constant* $M > 0$ *so that for all* $f \in H$ *we have:*

$$(6.1) \quad \sum_i |<f, f_i - g_i>|^2 \leq M \min\left(\sum_i |<f, f_i>|^2, \sum_i |<f, g_i>|^2\right).$$

*Moreover, if* $(g_i)$ *has a finite upper frame bound, then (1) and (2) are equivalent to*

*(3) there is a constant* $M > 0$ *so that for all* $f \in H$ *we have*

$$\sum_i |<f, f_i - g_i>|^2 \leq M \sum_i |<f, g_i>|^2.$$

*Proof.* $(1) \Rightarrow (2)$: Let $A_f$, $B_f$, $A_g$, $B_g$ be the lower and upper frame bounds of $(f_i)$ and $(g_i)$. Now, for all $f \in H$ we have

$$\sum_i |<f, f_i - g_i>|^2 \leq 2 \left( \sum_i |<f, f_i>|^2 + \sum_i |<f, g_i>|^2 \right) \leq$$

$$2 \left( \sum_i |f, f_i>|^2 + B_g \|f\|^2 \right) \leq 2 \left( \sum_i |<f, f_i>|^2 + \frac{B_g}{A_f} \sum_i |<f, f_i>|^2 \right)$$

$$\leq 2 \left( 1 + \frac{B_g}{A_f} \right) \sum_i |<f, f_i>|^2.$$

By symmetry, since we are assuming in (1) that $(g_i)$ is also a frame, we have

$$\sum_i |<f, f_i - g_i>|^2 \leq 2 \left( 1 + \frac{B_f}{A_g} \right) \sum_i |<f, g_i>|^2.$$

$(2) \Rightarrow (1)$: Given $M$ in (2) and any $f \in H$ we compute

$$A_f \|f\|^2 \leq \sum_i |<f, f_i>|^2 \leq 2 \left( \sum_i |<f, f_i - g_i>|^2 + \sum_i |<f, g_i>|^2 \right)$$

$$\leq 2 \left( M \sum_i |<f, g_i>|^2 + \sum_i |<f, g_i>|^2 \right) = 2(M+1) \sum_i |<f, g_i>|^2$$

$$\leq 2(M+1) \left( \sum_i |<f, f_i - g_i>|^2 + \sum_i |<f, f_i>|^2 \right) \leq$$

$$\leq 2(M+1) \left( M \sum_i |<f, f_i>|^2 + \sum_i |<f, f_i>|^2 \right)$$

$$\leq 2(M+1)^2 \sum_i |<f, f_i>|^2 \leq 2(M+1)^2 B_f \|f\|^2.$$

For the *moreover* part of the theorem, it is clear that (2) always implies (3). To see that (3) implies (1), we use the first four lines of the proof of $(2) \Rightarrow (1)$ to get a lower frame bound for $(g_i)$ and the upper frame bound comes from the fact that $(g_i)$ is Hilbertian. $\qquad \square$

There is also a Banach space version of Theorem 6.3.

**Theorem 6.4.** *Let $(f_i)$ and $(g_i)$ be sequences in a (Banach) Hilbert space $H$. The following are equivalent:*

*(1) $(f_i) \approx (g_i)$ (i.e. $Tf_i = g_i$ is a well-defined (possibly into) isomorphism of $H$ into $H$).*

*(2) There is a constant $M > 0$ so that for all sequences of scalars $(a_i)$ we have*

$$\| \sum_i a_i(f_i - g_i) \| \leq M \min \left( \| \sum_i a_i f_i \|, \| \sum_i a_i g_i \| \right).$$

*Moreover, if (2) holds, the equivalence constant in (1) is $\leq M + 1$.*

*Proof.* $(1) \Rightarrow (2)$: Given $Tf_i = g_i$ an isomorphism, for any finitely non-zero sequence of scalars $(a_i)$ let $f = \sum_i a_i f_i$. Then

$$\|f - Tf\| \leq \|f\| + \|Tf\| \leq (1 + \|T\|)\|f\|.$$

Similarly, $\|f - Tf\| \leq (1 + \|T^{-1}\|)\|Tx\|$. This is enough for (2).

$(2) \Rightarrow (1)$: By assumption, for any sequence of scalars $(a_i)$ we have

$$\| \sum_i a_i f_i \| \leq \| \sum_i a_i(f_i - g_i) \| + \| \sum_i a_i g_i \| \leq (M+1)\| \sum_i a_i g_i \|.$$

By symmetry we have that

$$\| \sum_i a_i g_i \| \leq (M+1)\| \sum_i a_i f_i \|.$$

$\square$

The disadvantage of Theorem 6.4 is that we do not get to conclude that the sequence $(g_i)$ is complete (i.e. spans $H$). For example, if $(e_i)$ is an orthonormal basis for $H$ and we let $f_i = e_i$ and $g_i = e_{i+1}$, the theorem holds.

Some of the first modern Hilbert space frame perturbation theorems come from the PhD thesis of Heil (see [49]). The motivation for studying abstract frame perturbations grew of out results of Feichtinger and Gröchenig and Walnut for perturbing WH-frames. This topic was then carried further by Christensen et al [32, 33, 45]. Good use of these results was made by Balan [8] in studying Fourier frames and Wavelet bases. We present a generalization due to Casazza and Christensen [32]. We will not track the frame bounds for $(g_i)$ below, but point out that in general it is quite important to know precisely the exact frame bounds obtained for a perturbation of a frame.

**Theorem 6.5.** *Let $(f_i)$ be a frame for a Hilbert space $H$ with frame bounds $A, B$. Let $(g_i)$ be a sequence in $H$ and assume $\max[\lambda_1 + \frac{\mu}{\sqrt{A}}, \lambda_2] < 1$. If one of the following two conditions is fulfilled for all finitely non-zero sequences of scalars $(a_i)$ and all $f \in H$, then $(g_i)$ is also a frame for $H$:*

$$\left( \sum_i |<f, f_i - g_i>|^2 \right)^{1/2} \leq \lambda_1 \left( \sum_i |<f, f_i>|^2 \right)^{1/2} + \mu\|f\|$$

$$\| \sum_i a_i(f_i - g_i) \| \leq \lambda_1 \| \sum_i a_i f_i \| + \lambda_2 \| \sum_i a_i g_i \| + \mu \sum_i |a_i|^2)^{1/2}.$$

*Proof.* For the proof, we refer the reader to the listed papers, or point out that these results can be obtained by judicious use of Theorems 6.3, 6.4. □

There are special perturbation theorems for WH-frames as well as exponential frames (see the article by Christensen in [74]). This ends our small introduction to perturbation theory. We mention that there is room for many useful perturbation results in this area. It would be particularly important, however, to produce good applications for new perturbation results so that this topic does not degenerate into total abstract nonsense. Good places to look for applications are WH-frames, frames of translates (See Section 8) and exponential frames (See Section 8).

## 7. WAVELET FRAMES

The theory of wavelet frames is one of the most underdeveloped areas of frame theory. (We should not confuse this with the study of orthogonal or biorthogonal wavelets, which is one of the most active research areas today.) Since these frames are also important for applications, it would be quite useful for much more progress to be made in this direction.

When the theory of wavelets was just blossoming, Daubechies, Grossman, and Meyer [57] combined the theory of the continuous wavelet transform with the theory of frames to define **affine frames** (or alternatively, **wavelet frames**) for $L^2(\mathbb{R})$. Daubechies [54] developed these ideas much further. As we have seen in the Balian-Low Theorem, it is not possible to have a Weyl-Heisenberg frame $(g, 1, 1)$ with $g$ both smooth and rapidly decaying. Some of the excitement over wavelets is due to a major result of Meyer [116] which says such phenomenon can occur there (and hence also in wavelet frames).

**Definition 7.1.** *Given* $g \in L^2(\mathbb{R})$, $a > 1$ *and* $b > 0$, *we say that* $(g, a, b)$ *generates an* **affine frame** *(or a* **wavelet frame***) for* $L^2(\mathbb{R})$ *if*

$$\{D_{a^n} T_{mb} g\}_{m,n \in \mathbb{Z}}$$

*is a frame for* $L^2(\mathbb{R})$. *The function* $g$ *is called the* **mother wavelet**. *The numbers* $a, b$ *are the* **frame parameters** *with* $a$ *being the* **dilation parameter** *and* $b$ *the* **shift parameter**.

It is sometimes necessary to use two mother wavelets in order to form a frame for $L^2(\mathbb{R})$. This was the original approach of [57].

**Theorem 7.2.** *Let* $g_1, g_2 \in L^2(\mathbb{R})$ *satisfy:*

(1) *supp* $(\hat{g}_1) \subset [-L, -\ell]$ *and supp* $(\hat{g}_2) \subset [\ell, L]$, *where* $0 < \ell < L < \infty$.

(2) $\hat{g}_1$ *and* $\hat{g}_2$ *are continuous and do not vanish on* $(-L, -\ell)$ *and* $(\ell, L)$, *respectively. Then* $(D_{a^n} T_{mb} g_1, D_{a^n} T_{mb} g_2)$ *is a frame for* $L^2(\mathbb{R})$ *for all* $1 < a < L/\ell$ *and all* $0 < b \leq 1/(L-1)$.

We really want to find only one mother wavelet generating a wavelet frame for $L^2(\mathbb{R})$. If the world were perfect, we could just add the two functions in Theorem 7.2 and get a good mother wavelet. Unfortunately, this does not work in general. For example, if $a = 2$, $b = 1$, $\hat{g}_1 = \chi_{(-2,-1]}$, and $\hat{g}_2 = \chi_{[1,2)}$ then $(D_{a^n}T_{mb}g_1, D_{a^n}T_{mb}g_2)$ actually forms an orthonormal basis for $L^2(\mathbb{R})$, but $\hat{f} = \chi_{(-2,-1]} - \chi_{[1,2)}$ is orthogonal to every function of the form $D_{a^n}T_{mb}(g_1 + g_2)$. The problem here is not critical however, since if we take $b$ small enough then $g_1 + g_2$ will generate a frame. This is an interesting result of Heil and Walnut [96].

**Theorem 7.3.** *Let $g_1, g_2 \in L^2(\mathbb{R})$ be as in Theorem 7.2. If $2L < 1/b$ then $(g_1 + g_2, a, b)$ generates an affine frame for $L^2(\mathbb{R})$.*

*Proof.* Let $g = g_1 + g_2$. Since supp $(D_{a^n}\hat{f} \cdot \hat{g}) \subset [\frac{-1}{2b}, \frac{1}{2b}]$, we have $D_{a^n}\hat{f} \cdot \hat{g} \in L^2[\frac{-1}{2b}, \frac{1}{2b}]$. Since $(b^{1/2}E_{mb})$ is an orthonormal basis for $L^2[\frac{-1}{2b}, \frac{1}{2b}]$, we have

$$\sum_{m,n} | < f, D_{a^n}T_{mb}g > |^2 = \sum_{m,n} | < D_{a^n}\hat{f} \cdot \overline{\hat{g}}, E_{-mb} > |^2$$

$$= b^{-1} \sum_n \int_R a^{-n}|\hat{f}(a^{-n}\gamma)|^2|\hat{g}(\gamma)|^2 d\gamma$$

$$= \int_{-\infty}^0 |\hat{f}(\gamma)|^2 \cdot b^{-1} \sum_n |\hat{g_1}(a^n\gamma)|^2 d\gamma + \int_0^\infty |\hat{f}(\gamma)|^2 \cdot b^{-1} \sum_n |\hat{g_2}|^2 d\gamma.$$

The result follows immediately from here. $\qquad\square$

Daubchies [54] has also given a general criterion for $g, a, b$ to generate a wavelet frame.

**Theorem 7.4.** *Let $g \in L^2(\mathbb{R})$ and $a > 1$ satisfy:*
*(1) There are positive constants $A, B$ such that*

$$A \le \sum_n |\hat{g}(a^n\gamma)|^2 \le B, \ a.e. \ \gamma \in \mathbb{R}.$$

*(2) $\lim_{b \to 0} \sum_{k \neq 0} \beta(k/b)^{1/2}\beta(-k/b)^{1/2} = 0$, where*

$$\beta(s) = ess \ sup_{|\gamma| \in [1,a]} \sum_n |\hat{g}(a^n\gamma)\hat{g}(a^n\gamma - s)|.$$

*Then there exists a $b_0 > 0$ such that $(g, a, b)$ generates a wavelet frame for $L^2(\mathbb{R})$ for each $0 < b < b_0$.*

Again, relying on the CC-condition, Casazza and Christensen [36] gave a stronger condition which also works for frame sequences (i.e. when $(g, a, b)$ generates a wavelet frame for its closed linear span).

**Theorem 7.5.** *Let $a > 1$, $b > 0$ and $g \in L^2(\mathbb{R})$ be given. Let*

$$N =: \{\gamma \in [1, a] : \sum_{n \in \mathbb{Z}} |\hat{g}(a^n \gamma)|^2 = 0\},$$

*and assume*

$$A =: inf_{|\gamma| \in [1,a]-N} \left[ \sum_{n \in \mathbb{Z}} |\hat{g}(a^n \gamma)|^2 - \sum_{k \neq 0} \sum_{n \in \mathbb{Z}} |\hat{g}(a^n \gamma)\hat{g}(a^n \gamma + k/b)| \right] > 0,$$

$$B =: sup_{|\gamma| \in [0,a]} \sum_{k,n \in \mathbb{Z}} |\hat{g}(a^n \gamma)\hat{g}(a^n \gamma + k/b)| < \infty.$$

*Then $\{\frac{1}{a^{n/2}} g(\frac{x}{a^n} - mb)\}_{n,m \in \mathbb{Z}}$ is a frame sequence with bounds $\frac{A}{b}, \frac{B}{b}$.*

As we saw earlier, the dual of a WH-frame is another WH-frame. One major difficulty with wavelet frames is that the dual frame of a wavelet frame need not be a wavelet frame. What we do have is that given a wavelet frame with frame operator $S$, $S(D_{a^n} f) = D_{a^n} S f$ and similarly for $S^{-1}$. An important problem, especially for applications, is to classify the affine frames which have duals which are also affine frames. Recently Bownik [26] gave two equations which are necessary and sufficient for two affine frames to be dual frames. We ask:

**Problem 7.6.** *Classify all $(g, a, b)$ which generate wavelet frames.*

The above problem may be completely intractable at this time. Something which would be just as good in practice is:

**Problem 7.7.** *Find all wavelet frames whose cannonical dual frame is also a wavelet frame (and preferably give an exact representation for the dual).*

As we saw in Section 5, a WH-frame is a Riesz basis for $L^2(\mathbb{R})$ if and only if $ab = 1$. Moreover, in this case, functions giving WH-frames are either not smooth or do not decay quickly. Meyer [116] showed that things are much better in the wavelet frame case by exhibiting a $C^\infty$ function with compactly supported Fourier transform which generates a wavelet basis for $L^2(\mathbb{R})$.

**Definition 7.8.** *(The **Meyer wavelet**) The **Meyer wavelet** is the function $\psi \in L^2(\mathbb{R})$ defined by $\hat{\psi}(\gamma) = e^{i\gamma/2}\omega(|\gamma|)$, where*

$$\omega(\gamma) = \begin{cases} 0 : & \gamma \leq \frac{1}{3} \text{ or } \gamma \geq \frac{4}{3} \\ sin \frac{\pi}{2}\nu(3\gamma - 1) : & \frac{1}{3} \leq \gamma \leq \frac{2}{3} \\ cos \frac{\pi}{2}\nu(\frac{3\gamma}{2} - 1) : & \frac{2}{3} \leq \gamma \leq \frac{4}{3} \end{cases}$$

*and $\nu \in C^\infty(\mathbb{R})$ is such that $\nu(\gamma) = 0$ for $\gamma \leq 0$, $\nu(\gamma) = 1$ for $\gamma \geq 1$, $0 \leq \nu(\gamma) \leq 1$ for $\gamma \in [0, 1]$, and $\nu(\gamma) + \nu(1 - \gamma) = 1$ for $\gamma \in [0, 1]$.*

It takes some work to show that $(D_{2^n}T_m\psi)_{m,n\in\mathbb{Z}}$ is an orthonormal basis for $L^2(\mathbb{R})$, and any $\psi$ with the property that $(D_{2^n}T_m\psi)$ forms an orthonormal basis for $L^2(\mathbb{R})$ is called an **(orthonormal) wavelet**. Mallat and Meyer then developed the theory of multiresolution analysis to put such examples into a natural framework. This was based on some earlier remarkable results of Mallat [115].

**Definition 7.9.** *A* **multiresolution analysis** *for $L^2(\mathbb{R})$ consists of*
 (1) *Closed subspaces $V_n \subset L^2(\mathbb{R})$ for $n \in \mathbb{Z}$ and satisfying:*
    *(a)* $V_n \supset V_{n+1}$,
    *(b)* $\cap V_n = \{0\}$,
    *(c)* $\cup V_n$ *is dense in $L^2(\mathbb{R})$,*
    *(d)* $V_{n+1} = D_2 V_n = \{D_2 f : f \in V_n\}$,
 (2) *A function $\phi \in V_0$ such that $(T_m\phi)_{m\in\mathbb{Z}}$ is an orthonormal basis for $V_0$.*

The important point is that each MRA generates a natural orthonormal basis for $L^2(\mathbb{R})$. There is an MRA generating the Meyer wavelet. The Haar system is generated by setting

$V_0 = \{f \in L^2(\mathbb{R}) : \text{f is a constant on each of the intervals } [m, m+1), \ m \in \mathbb{Z}\}$

and letting $V_n = D_{2^n}V_0$, and $\phi = \chi_{[0,1)}$.

How does an MRA work? Since $V_n \subset V_{n-1}$, there is a subspace $W_n$ of $V_{n-1}$ so that $V_{n-1} = V_n \oplus W_n$. By assumption there is a function $\phi \in W_0$ such that $(T_m\phi)$ is an orthonormal basis for $W_0$. Now, $W_{n+1} = D_2 W_n$ making $(D_{2^n}T_m\phi)$ an orthonormal basis for $W_n$. Since $L^2(\mathbb{R}) = (\sum \oplus W_n)_{\ell_2}$, it follows that $(D_{2^n}T_m\phi)$ is an orthonormal basis for $L^2(\mathbb{R})$.

Multiresolution analysis is one of the most used and useful tools in wavelet theory. We refer the reader to the classical book of Daubechies [55] for an introduction to this important topic. Also, more recently, the notion of a **frame multiresolution analysis** has been defined by Benedetto and Li [18]. This topic is just now being developed and we refer the reader to [18, 20, 105] for the most recent results in this area.

Once again we emphasize that the focus of these notes is on modern abstract frame theory. This accounts for the sparse treatment of wavelets and wavelet frames. The reader should be aware however that this is one of the largest and most active areas of research today.

## 8. Frames of Translates

Frames consisting of translates of a single function play an important role in sampling theory as well as in wavelet theory and Weyl-Heisenberg frame theory. In this section we will be dealing with families $(f_i)$ in $L^2(\mathbb{R})$ which are not frames for $L^2(\mathbb{R})$, but are frames for their closed linear span. We will call such sequences **frame sequences**. It is known that no sequence of translates of a single function can form a frame for $L^2(\mathbb{R})$. This is a recent result of Christensen, Deng and Heil [48]. Previously, Olson and Zalik [120] showed there are no Riesz bases of translates for $L^2(\mathbb{R})$.

We first introduce some notation. If $\phi \in L^2(\mathbb{R})$ and $b > 0$, we define the function $\Phi_b : \mathbb{T} \to \mathbb{R}$ by

$$\Phi_b(\zeta) = \sum_{n \in \mathbb{Z}} |\hat{\phi}\left(\frac{\zeta + n}{b}\right)|^2.$$

Note that $\Phi_b \in L^1(\mathbb{T})$. A direct calculation shows that for any $n \in \mathbb{Z}$ we have

$$< T_{nb}\phi, \phi > = < e^{-2\pi i n \zeta b}\hat{\phi}, \hat{\phi} > = \frac{1}{b}\int_0^1 \Phi_b(\zeta)e^{-2\pi i n \zeta}d\zeta = \frac{1}{b}\hat{\Phi}_b(n).$$

If $\Lambda \subset \mathbb{Z}$, we let $H_\Lambda$ be the closed subspace of $L^2(\mathbb{T})$ generated by the characters $e^{2\pi i n \zeta}$, for $n \in \Lambda$. We let $E_\Lambda$ be the closed subspace of $H_\Lambda$ consisting of all $f$ such that $\Phi_b(\zeta)f(\zeta) = 0$ a.e. If $f \in H_\Lambda$, we denote by $d(f, E_\Lambda)$ the distance of $f$ to the subspace $E_\Lambda$. Our first result is due to Casazza, Christensen and Kalton [40]

**Theorem 8.1.** *Suppose $\phi \in L^2(\mathbb{R})$ and $b > 0$. If $\Lambda \subset \mathbb{Z}$ then $(T_{nb}\phi)_{n \in \Lambda}$ is a frame sequence with frame bounds $A$ and $B$ if and only if for every $f \in H_\Lambda$ we have*

$$Ad(f, E_\Lambda)^2 \leq \frac{1}{b}\int_0^1 |f(\zeta)|^2\Phi_b(\zeta)d\zeta \leq B\|f\|^2,$$

*or equivalently, for all $f \in H_\Lambda \cap E_\Lambda^\perp$,*

$$A\|f\|^2 \leq \frac{1}{b}\int_0^1 |f(\zeta)|^2\Phi_b(\zeta)d\zeta \leq B\|f\|^2.$$

*Furthermore, if this condition is satisfied, then $(T_{nb}\phi)_{n \in \Lambda}$ is an exact frame sequence with the same frame bounds if and only if $E_\Lambda = \{0\}$.*

*Proof.* By our earlier remarks and the definitions, $(T_{nb}\phi)_{n \in \Lambda}$ is a frame sequence with frame bounds $A, B$ if and only if the linear map $T : c_{00}(\Lambda) \to L^2(\mathbb{R})$ (where $c_{00}(\Lambda)$ denotes the finitely non-zero sequences supported on $\Lambda$) defined by $Te_n = T_{nb}\phi$ extends to a bounded linear operator $T : \ell^2 \to L^2(\mathbb{R})$ with

$$Ad(u, \ker T)^2 \leq \|Tu\|^2 \leq B\|u\|^2,$$

for all $u \in \ell^2(\Lambda)$.

Let $U : H_\Lambda \to \ell^2(\Lambda)$ be the natural isometry $Uf = \{\hat{f}(n)\}_{n \in \Lambda}$. Then for any trigonometric polynomial $f \in H_\Lambda$ we have

$$\|TUf\|^2 = \|\sum_{n \in \Lambda} \hat{f}(n)T_{nb}\phi\|^2 = \int_{-\infty}^\infty |\sum_{n \in \Lambda} \hat{f}(n)e^{-2\pi i n b \zeta}\hat{\phi}(\zeta)|^2 d\zeta$$

$$= \frac{1}{b} \sum_{n \in \mathbb{Z}} \int_0^1 |f(\zeta)|^2 |\hat{\phi}(\frac{n+\zeta}{b})|^2 d\zeta = \frac{1}{b} \int_0^1 |f(\zeta)|^2 \Phi_b(\zeta) d\zeta.$$

The theorem now follows. $\qquad \square$

Theorem 8.1 yields a generalization of a result of Benedetto and Li [18]. The original theorem had an unnecessary hypothesis which our theorem removes. However, Kim and Lim [105] also removed this unnecessary hypothesis as did Benedetto and Treiber [19].

**Theorem 8.2.** *If $\phi \in L^2(\mathbb{R})$ and $b > 0$ then*

*(1) $(T_{nb}\phi)_{n \in \mathbb{Z}}$ is an orthonormal sequence if and only if*

$$\Phi_b(\gamma) = b \ \ a.e.$$

*(2) $(T_{nb}\phi)_{n \in \mathbb{Z}}$ is an exact frame sequence with frame bounds $A, B$ if and only if*

$$bA \leq \Phi_b(\gamma) \leq bB \ \ a.e.$$

*(3) $(T_{nb}\phi)_{n \in \mathbb{Z}}$ is a frame sequence with frame bounds $A, B$ if and only if*

$$bA \leq \Phi_b(\gamma) \leq bB \ \ a.e.$$

*on $\mathbb{T} - N_b$ where $N_b = \{\zeta \in \mathbb{T} : \Phi_b(\zeta) = 0\}$.*

*Proof.* Note that (1) follows easily from the fact that $(T_{nb}\phi)_{n \in \mathbb{Z}}$ is orthonormal if and only if $TU$ is unitary. (2) is immediate from (3). For (3) we note that if $\Lambda = \mathbb{Z}$, then $H_\Lambda = L^2(\mathbb{T})$ and $E_\Lambda = L^2(N_b)$. Hence, $d(f, E_\Lambda)^2) = \int_{\mathbb{T} - N_b} |f|^2 d\zeta$. This is all we need. $\qquad \square$

We now state an interesting consequence of the above ideas.

**Theorem 8.3.** *Suppose $\phi \in L^2(\mathbb{R})$ and $b > 0$. Then the sequence $(T_{nb}\phi)_{n \in \mathbb{N}}$ is a frame sequence if and only if $(T_{nb}\phi)_{n \in \mathbb{N}}$ is a Riesz basis for its span.*

*Proof.* Assume that $(T_{nb}\phi)_{n \in \mathbb{N}}$ is a frame sequence. then (see for example [64]) if $0 \neq f \in H_\mathbb{N}$ we have that $\log |f| \in L^1$ and so, in particular, $|f| > 0$ a.e. This implies that $E_\mathbb{N} = 0$. It follows that if $A, B$ are the frame bounds for $(T_{nb}\phi)_{n \in \mathbb{N}}$ then for every trigonometric polynomial in $H_\mathbb{N}$ we have

$$A\|f\|^2 \leq \int |f(\zeta)|^2 \Phi_b(\zeta) d\zeta \leq B\|f\|^2.$$

Now, if $f$ is any trigonometric polynomial, then for large enough $n$ we have that $e^{2\pi i n \zeta} f \in H_\mathbb{N}$. Thus, the same inequality follows for all trigonometric polynomials in $L^2(\mathbb{T})$. This implies the theorem. $\qquad \square$

A careful examination of the proof of Theorem 8.3 should convince the reader that the theorem holds for $(T_{nb}\phi)_{n \in \Lambda}$ with $\Lambda \subset \mathbb{N}$. Theorem 8.2 is a classification theorem for when an evenly spaced sequence of translates forms
a frame, Riesz basis or orthonormal basis. The main question here is:

**Problem 8.4.** *Classify those functions $\phi \in L^2(\mathbb{R})$ and sequences $\lambda_n \in \mathbb{R}$ so
that $(T_{\lambda_n}\phi)$ is a frame sequence.*

Problem 8.4 is way beyond the scope of the theory at this time. But, we can
break this problem down into smaller parts. First we ask whether every frame
of translates is equivalent to a subset of an evenly spaced frame of translates?

**Problem 8.5.** *Is every frame of translates equivalent to a frame sequence
$(T_{nb}\phi)_{n \in \Lambda}$ for some $\Lambda \subset \mathbb{Z}$ (perhaps with multiples of the elements)?*

It can be shown that the answer to Problem 8.5 is *yes* if $\hat{\phi}$ is bounded. This
is a simple variation of the proof of Theorem 5.3.2 in the section of Christensen
of the book [74]. Problem 8.5 could benefit from some deeper perturbation
theory.

**Problem 8.6.** *Classify the frame sequences consisting of subsets of evenly
spaced frames of translates.*

We end with a few remarks about **exponential frames**. Here we are in-
terested in finding the $\gamma$ and families of real (or complex) numbers $\lambda_n$ so that
$(e^{i\lambda_n})$ forms a frame (or a Riesz basis) for $L^2(-\gamma, \gamma)$. Since the Fourier trans-
form takes $T_{nb}\phi$ to $E_{-nb}\hat{\phi}$, exponential frames (or **Fourier frames** as they are
sometimes called) are related to frames of translates. This topic is so large
that entire books are devoted to it (see Young [139]). Therefore, we have
decided not to cover it at all since we cannot do justice to it in the limited
space we have here. There is a rich history to these questions going back to
Paley and Wiener (See for example Young [139]). Also, there is a large body
of information available concerning exponential frames and we refer the reader
to Young [139] for this or the bibliography of Seip [130] or Balan [8] for some
more recent results. Young's book also contains extensive stability results (i.e.
Perturbation results).

Another major direction of research related to this topic is sampling theory.
This has been an important topic since **Shannon's sampling theorem** first
appeared [132, 136]. This also involves a fundamental question in signal pro-
cessing. Namely, how do we represent a function on $\mathbb{R}$ in terms of a discrete
sequence? One way is to sample the function - either on a uniform grid or
an irregular grid. The Shannon-Whittaker Representation Theorem [132, 136]
says that the class of bandlimited functions (described in the theorem) can be
completely characterized by their sample values. Again, sampling theory is a
large area of research and we have chosen not to enter this important topic
here. We refer the reader to [2, 3, 14, 16, 21, 113] and their bibliographies for
an up to date view of this important topic.

## 9. Recent Developments in Weyl-Heisenberg Frame Theory

In this section we will look at some recent results in Weyl-Heisenberg frame theory. We will treat this section more as a survey of the latest results here.

An important result was proved recently and independently by Daubechies, H. Landau and Z. Landau [60], Janssen [102], and Ron and Shen [126].

**Theorem 9.1.** *For $g \in L^2(\mathbb{R})$ and $a, b > 0$, the following are equivalent:*
*(1) $(g, a, b)$ is a Weyl-Heisenberg frame.*
*(2) $(E_{m/a}T_{n/b}g)_{n,m\in\mathbb{Z}}$ is a Riesz basic sequence.*

Casazza and Lammers [42] split the above theorem in half and showed,

**Theorem 9.2.** *For $g \in L^2(\mathbb{R})$ and $a, b > 0$, the following are equivalent:*
*(1) There are $A, B > 0$ so that $A \leq \sum_{n\in\mathbb{Z}} |g(t - na)|^2 \leq B$, a.e.*
*(2) $(E_{m/a}g)_{m\in\mathbb{Z}}$ is a Riesz basic sequence.*

Next, we have the classification theorem for tight WH-frames. The equivalence of (1) and (2) in this theorem can be found in [126], Corollary 2.19, while it also follows from the developments in [103], Subsection 1.3. The equivalence of (1) and (3) can be deduced from [126], Corollary 6.8 or from [103], Theorem 1.4.1. Our proof comes from [38], Theorem 3.2.

**Theorem 9.3.** *Let $g \in L^2(\mathbb{R})$ and $a, b \in \mathbb{R}$. The following are equivalent:*
*(1) $(E_{mb}T_{na}g)$ is a normalized tight Weyl-Heisenberg frame for $L^2(\mathbb{R})$.*
*(2) We have:*
    *(a) $G_0(t) = \sum_{n\in\mathbb{Z}} |g(t - na)|^2 = b$ a.e.*
    *(b) For all $k \neq 0$, $G_k(t) = \sum_n g(t - na)\overline{g(t - na - k/b)} = 0$ a.e.*
*(3) We have $g \perp E_{n/a}T_{m/b}g$, for all $(n, m) \neq (0, 0)$ and $\|g\|^2 = ab$.*
*(4) $(E_{n/a}T_{m/b}g)$ is an orthogonal sequence in $L^2(\mathbb{R})$ and $\|g\|^2 = ab$.*
*(5) $(E_{mb}T_{na}g)$ is a Weyl-Heisenberg frame for $L^2(\mathbb{R})$ with frame operator $S$ and $Sg = g$.*

*Moreover, when at least one of $(1) - (5)$ holds, $(E_{mb}T_{na}g)$ is an orthonormal basis for $L^2(\mathbb{R})$ if and only if $\|g\| = 1$.*

*Proof.* (1) $\Leftrightarrow$ (2). Assume $(E_{mb}T_{na}g)$ is a normalized tight frame. For any function $f \in L^2(\mathbb{R})$ which is bounded and supported on an interval of length $< 1/b$, we have by the WH-frame Identity that $F_2(f) = 0$. Hence, again by the WH-frame identity we have

$$\|f\|^2 = \int_R |f(t)|^2 dt = \sum_{n,m} | < f, E_{mb}T_{na}g > |^2 = b^{-1} \int_R |f(t)|^2 G_0(t) dt.$$

(2) follows easily from here. The converse is immediate from the WH-frame Identity.

(3) $\Leftrightarrow$ (4). This is a direct calculation.

$(1) \Leftrightarrow (5)$. Since $S$ commutes with translation by $a$ and modulation by $b$, (5) is equivalent to $S = I$ which, in turn, is equivalent to (1).

$(2) \Leftrightarrow (3)$. If we compute:

$$< g, E_{mb}T_{na}g >= \int_R g(t)\overline{E_{mb}g(t-na)}dt = \int_R g(t)\overline{g(t-na)}e^{-2\pi imbt}dt =$$

$$\int_0^{1/b} \sum_{k\in\mathbb{Z}} g(t-k/b)\overline{g(t-na-k/b)}e^{-2\pi imbt}dt.$$

A little reflection should convince the reader that this is all we need.    $\square$

There are several consequences of Theorem 9.3.

**Remark 9.4.** *Another application of Theorem 9.3 (see [38], Corollary 3.4) is a simple proof of Theorem 5.19 (2).*

*Proof.* If $(E_{mb}T_{na}g)$ is a WH-frame with frame operator $S$, then by Proposition 5.5 and the remarks following equation (4.3), $(E_{mb}T_{na}S^{-1/2}g)$ is a normalized tight WH-frame equivalent to our WH-frame. So if our frame is exact then $(E_{mb}T_{na}S^{-1/2}g)$ is an orthonormal basis. Hence,

$$1 = \|S^{-1/2}g\|^2 = ab.$$

Also, if $ab = 1$, then $(E_{mb}T_{na}S^{-1/2}g)$ is a tight WH-frame and hence by Theorem 9.3, $(E_{n/a}T_{m/b}S^{-1/2}g)$ is an orthogonal sequence. But, $ab = 1$ implies that $n/a = nb$ and $m/b = ma$, so $(E_{nb}T_{ma}S^{-1/2}g)$ is an orthogonal basis for $L^2(\mathbb{R})$. Since $S^{-1/2}$ is an invertible operator, it follows that $(E_{mb}T_{na}g)$ is a Riesz basis.    $\square$

We are now in a position to write down explicitly the functions which yield WH-frames which are also orthonormal bases for $L^2(\mathbb{R})$. If $f(x,y)$ is a function of two variables, we denote by $f_x(y)$ (respectively, $f_y(x)$) the function of one variable given by $f_x(y) = f(x,y)$ (respectively, $f_y(x) = f(x,y)$).

**Theorem 9.5.** *Let $ab = 1$ and $g \in L^2(\mathbb{R})$. The following are equivalent:*

*(1) $(E_{mb}T_{na}g)$ is a normalized tight WH-frame for $L^2(\mathbb{R})$.*

*(2) $(E_{mb}T_{na}g)$ is an orthonormal basis for $L^2(\mathbb{R})$.*

*(3) There is a measurable function $f : [0,1] \times [o,a) \to \mathbb{R}$ such that if*

$$h(x,y) = \sqrt{b}e^{-2\pi i f(x,y)}$$

*then for all $n \in \mathbb{Z}$ and all $y \in [0,a)$ we have*

$$g(y+na) = \widehat{h_y}(n).$$

*Proof.* Since $ab = 1$, we have $a = 1/b$, so Theorem 9.3 (2) becomes:

(a) $G(t) = \sum_{n \in \mathbb{Z}} |g(t - na)|^2 = b$ a.e.

(b) $G_k(t) = \sum_{n \in \mathbb{Z}} g(t - na)\overline{g(t - (n-k)a)} = 0$, a.e. for all $k \neq 0$.

But, condition (b) is equivalent to $z_y = (g(y - na))_{n \in \mathbb{Z}}$ is orthogonal to all of its proper shifts and (a) is equivalent to $\|z_y\|^2 = b$. Direct calculation (see the proof of Theorem 4.1, [38]) yields: For each $y \in [0, a)$ there is a function $f_y : [0, 1] \to \mathbb{R}$ and $h_y : [0, a) \to C$ with

$$h_y(x) = e^{-2\pi i f_y(x)},$$

and

$$\widehat{h_y}(n) = g(y + na).$$

So defining $f(x, y) : [0, 1] \times [0, a) \to \mathbb{R}$ and $h(x, y)$ by:

$$f(x, y) = f_y(x), \qquad h(x, y) = h_y(x),$$

yields our result modulo the measurability conditions which are obvious. □

A little reflection should convince the reader that the function $h(x, y)$ above is essentially just the Zak transform in disguise. We can also quickly get from here most of the standard examples of WH-frames from the literature as well as some new ones. We are working with $ab = 1$ in these examples.

**Example 9.6.** *Letting $f(x, y) = k$, for all $y \in [0, a)$ gives the function*

$$g(x) = \sqrt{b}\chi_{[0,a)}.$$

*Now, $(g, a, b)$ generates an orthonormal basis for $L^2(\mathbb{R})$.*

**Example 9.7.** *Letting*

$$f(x, y) = K_y, \quad \text{for all} \quad x \in [0, 1],$$

*gives the function*

$$g(x) = \sqrt{b}e^{ih(x)}\chi_{[0,a)},$$

*where $h : [0, 1] \to \mathbb{R}$ yielding that $(g, a, b)$ generates an orthonormal basis for $L^2(\mathbb{R})$.*

**Example 9.8.** *Letting $f(x, y) = f(x, y')$, for all $y, y' \in [0, a]$ gives the function*

$$g(x) = \sum_{n \in \mathbb{Z}} c_n \chi_{[na,(n+1)a)},$$

*where $(c_n)$ is an $\ell_2$-sequence which is orthogonal to all its proper shifts. So $(g, a, b)$ generates an orthonormal basis for $L^2(\mathbb{R})$.*

**Example 9.9.** *If we partition $[0, a)$ into disjoint measurable sets $(A_n)_{n \in \mathbb{Z}}$, let $B_n = A_n + \{n\}$ and define*

$$f(x, y) = d_n + nx, \quad \text{for all} \quad y \in A_n,$$

*where $d_n \in C$, we get*

$$g(x) = \sum_{n \in \mathbb{Z}} c_n \chi_{B_n}.$$

*Thus, $(g, a, b)$ generates an orthonormal basis for $L^2(\mathbb{R})$.*

We can generalize Theorem 9.5 to the case $1/ab \in \mathbb{Z}$ [38] and perhaps even to arbitrary rational $ab$. However, we ask

**Problem 9.10.** *Give an explicit representation of all functions $g \in L^2(\mathbb{R})$ and $ab$ irrational so that $(g, a, b)$ generates a normalized tight WH-frame.*

All normalized tight Weyl-Heisenberg frames have the same frame operator - namely, the identity operator. It is possible in general to classify the WH-frames with the same frame operator. As we will see, this result is a natural generalization of the equivalence of (1) and (2) in Theorem 9.3. Since we now have to work with two WH-frames, we define for $(g, a, b)$ and $k \in \mathbb{Z}$

$$(9.1) \qquad G_k^g = \sum_{n \in \mathbb{Z}} g(t - na) \overline{g(t - na - k/b)}.$$

Our next theorem comes from [39].

**Theorem 9.11.** *The systems $(g, a, b)$ and $(h, a, b)$ generate WH-frames with the same frame operator if and only if $G_k^g = G_k^h$ a.e. for all $k \in \mathbb{Z}$.*

In [39], Theorem 8.2, there is a generalization of Theorem 9.11 to the case where we have $(g, a, b)$ and $(h, c, d)$.

As we saw in Proposition 5.10, if $(E_{mb}T_{na}g)$ is a WH-frame, then the frame operator $S$ is given by

$$Sf = \sum_{n,m \in \mathbb{Z}} <f, E_{mb}T_{na}g> E_{mb}T_{na}g.$$

We know that the frame operator is an invertible operator on $L^2(\mathbb{R})$ and that both $S$ and $S^{-1}$ commute with translation by $a$ and modulation by $b$. Walnut [134] gave a useful explicit representation for $S$ even for many cases where this is only a preframe operator.

**Theorem 9.12.** *Let $a, b > 0$ and $g \in W(L^\infty, L^1)$. Then the sum $Sf$ converges unconditionally for each $f \in L^2(\mathbb{R})$ and is given by*

$$(9.2) \qquad Sf = b^{-1} \sum_k T_{k/b}f \cdot G_k.$$

The sum in (9.2) is called the **Walnut representation** or the **Walnut series** for the frame operator and it also converges unconditionally. A detailed study of the convergence properties of the Walnut representation of the frame operator was done by Casazza, Christensen and Janssen [39]. Recall from Proposition 5.33, for $M \subset \mathbb{Z}$ we define the operator $S_M : L^2(\mathbb{R}) \to L^2(\mathbb{R})$ by

$$S_M f = \sum_{k \in M} f(\cdot - k) G_k = \sum_{k \in M} (T_k f) \cdot G_k.$$

If $L, K \in \mathbb{Z}^+$ and $M = \{-L, -L+1, \cdots, 0, 1, \cdots K\}$ we write $S_{L,K} = S_M$. Also, we write $S_K =: S_{K,K}$. We say that the Walnut representation converges **weakly** (respectively, **in norm**) if for every $f \in L^2(\mathbb{R})$, $S_{L,K} f \to S f$ weakly (respectively, in norm) as $K, L \to \infty$. If $S_M \to Sf$ as a net for every $f \in L^2(\mathbb{R})$ we say the Walnut series converges **unconditionally**. Finally, for each of the above, we can discuss the **symmetric** convergence of the Walnut representation. i.e. If $S_K f \to Sf$ for all $f \in L^2(\mathbb{R})$.

We saw earlier that there are WH-frames for which the CC-condition fails. Also [39] there are WH-frames for which the Walnut representation does not converge even weakly for some $f \in L^2(\mathbb{R})$. However, weak and norm (strong) convergence of the Walnut representation are the same [39] as we now show.

**Theorem 9.13.** *For $a, b \in \mathbb{R}$ and $g \in L^2(\mathbb{R})$ with $|G_o(t)| \leq B$ a.e., the following are equivalent:*

*(1) The Walnut series converges in norm (resp. symmetrically in norm) for every $f \in L^2(\mathbb{R})$.*

*(2) the Walnut series converges weakly (resp. weakly symmetrically) for every $f \in L^2(\mathbb{R})$.*

*(3) We have $\sup_{L,K} \|S_{L,K}\| = B < \infty$ (resp. $\sup_K \|S_K\| = B < \infty$).*

*Proof.* Corollary 5.13 says that the Walnut series converges unconditionally on a dense subset of $L^2(\mathbb{R})$. Therefore, weak and norm convergence become equivalent to the partial sum operators being uniformly bounded. □

Symmetric and norm convergence of the Walnut representation are not equivalent in general [39]. It is natural to ask:

**Problem 9.14.** *Give necessary and sufficient conditions, in terms of $g, a, b$ for the Walnut series to converge (symmetrically) in norm for all $f \in L^2(\mathbb{R})$.*

It is known that unconditional convergence of the Walnut representation is not equivalent to norm convergence (see [39]). We do have the corresponding result to Theorem 9.13 for unconditional convergence [39]. The proof is similar to that case combined with Theorem 2.3.

**Theorem 9.15.** *Let $a, b \in \mathbb{R}$ and $g \in \mathbf{PF}$. The following are equivalent:*

*(1) The Walnut series converges weakly unconditionally for every $f \in L^2(\mathbb{R})$.*

(2) *the Walnut series converges unconditionally in norm for every* $f \in L^2(\mathbb{R})$.

(3) $\sup_{M \subset Z, |M| < \infty} \|S_M\| < \infty$.

It is also known that the CC-condition is strong enough to yield the unconditional convergence of the Walnut series. The following is a compilation of several results from [39].

**Theorem 9.16.** *Let* $a, b \in \mathbb{R}$ *and* $g \in \mathbf{PF}$.

(1) *If* $(g, a, b)$ *satisfies the CC-condition then the Walnut series converges unconditionally for all* $f \in L^2(\mathbb{R})$.

(2) *If* $ab$ *is rational, then the following are equivalent:*

  (i) *The Walnut series converges unconditionally for all* $f \in L^2(\mathbb{R})$.

  (ii) *There is a* $B > 0$ *so that* $\sum_k |G_k(t)| \le B$ *a.e.*

(3) *If* $ab$ *is irrational, then there is a WH-frame* $(g, a, b)$ *for which* (i) *and* (ii) *in* (2) *are not equivalent.*

*Proof.* We will prove (1) to give some of the flavor of this result. The other parts are more technical and we refer the reader to [39] for these. For any $h \in L^2(\mathbb{R})$ and any $m \ge n > 0$ we have (using the fact that $T_{-k/b} G_k = \overline{G_{-k}}$)

$$(9.3) \qquad \sum_{|k|=n}^{m} |<(T_{k/b}f)G_k, h>| = \sum_{|k|=n}^{m} |\int \overline{h(t)}(T_{k/b}f)(t)G_k(t)dt| \le$$

$$\sum_{|k|=n}^{m} \int |h(t)||T_{k/b}f(t)||G_k(t)|dt \le$$

$$\sum_{|k|=n}^{m} \left( \int |h(t)|^2 |G_k(t)|dt \right)^{1/2} \cdot \left( \int |T_{k/b}f(t)|^2 |G_k(t)|dt \right)^{1/2} \le$$

$$\left( \sum_{|k|=n}^{m} \int |h(t)|^2 |G_k(t)|dt \right)^{1/2} \cdot \left( \sum_{|k|=n}^{m} \int |T_{k/b}f(t)|^2 |G_k(t)|dt \right)^{1/2} =$$

$$\left( \int |h(t)|^2 \sum_{|k|=n}^{m} |G_k(t)|dt \right)^{1/2} \cdot \left( \int |f(t)|^2 \sum_{|k|=n}^{m} |T_{-k/b}G_k(t)|dt \right)^{1/2} =$$

$$\left( \int |h(t)|^2 \sum_{|k|=n}^{m} |G_k(t)|dt \right)^{1/2} \cdot \left( \int |f(t)|^2 \sum_{|k|=n}^{m} |G_{-k}(t)|dt \right)^{1/2}.$$

By our hypotheses,

$$|h|^2 \sum_k |G_k| \in L^1(R).$$

Hence, by the Lebesgue Dominated Convergence Theorem, we conclude that the following series converges in $L^1(R)$:

$$\sum_k |h|^2 |G_k|.$$

It follows that the right hand side of (9.3) goes to zero as $n \to \infty$. We conclude that the Walnut series for $f$ is weakly unconditionally Cauchy in $L^2(\mathbb{R})$ and hence it is unconditionally convergent in norm by the Theorem 2.3 □

Theorem 9.16 gives a classification of the WH-frames $(g, a, b)$ for which the Walnut series always converges unconditionally for the case $ab$ rational.

**Problem 9.17.** *Find necessary and sufficient conditions on $g, a, b$ with $ab$ irrational which guarantee that the Walnut series converges unconditionally for all $f \in L^2(\mathbb{R})$.*

It would also be interesting to know, for a given WH-system $(g, a, b)$ for which the Walnut representation does not always converge unconditionally, a large class of functions for which the Walnut representation does converge unconditionally anyway. One such class is discussed in Theorem 6.10 of [39].

If $(g, a, b)$ is a normalized tight WH-frame then the frame operator is $S = I$ and hence the frame operator extends to be a bounded linear operator from $L^p(\mathbb{R})$ to $L^p(\mathbb{R})$, for all $1 \le p \le \infty$. It is natural then to consider general conditions which allow the frame operator to extend to be a bounded linear operator (or an invertible operator) on $L^p(\mathbb{R})$. These $(g, a, b)$ were classified by Casazza, Christensen and Janssen ([39], Theorem 7.1).

**Theorem 9.18.** *If $ab \le 1$ and $g \in L^2(\mathbb{R})$, the following are equivalent:*
*(1) There is a constant $B > 0$ so that*

$$\sum_{k \in \mathbb{Z}} |G_k(t)| \le B, \quad a.e. \ x \in \mathbb{R}.$$

*(2) The frame operator $Sf = \sum_{n,m \in \mathbb{Z}} < f, E_{mb}T_{na}g > E_{mb}T_{na}g$ extends to a bounded linear operator from $L^p(\mathbb{R})$ to $L^p(\mathbb{R})$ for every $1 \le p \le \infty$. (Here, by $L^\infty$ we really mean $c_0(\mathbb{R})$ - the closure of the compactly supported functions in $L^\infty$.)*

*Moreover, if $(g, a, b)$ generates a WH-frame, $ab$ is rational, and $g$ satisfies the uniform CC-condition, then $S$ extends to an invertible operator on $L^p(\mathbb{R})$ for all $1 \le p \le \infty$.*

*Proof.* $(1) \Rightarrow (2)$ : We will show that $S$ is a bounded linear operator mapping a dense subset of $L^r(\mathbb{R})$ into itself for $r = 1, \infty$. Then, by the Riesz-Thorin Interpolation Theorem [141], page 95, $S$ extends to a bounded linear operator on $L^p(\mathbb{R})$ for all $1 \le p \le \infty$.

*Case I: The $L^1$-Case.*

If $f \in L^1$ is bounded and compactly supported, then by Corollary 5.13 we have:

$$\|Sf\|_{L^1} = \|Lf\|_{L^1} = b^{-1} \int_R |\sum_k (T_{k/b}f)G_k| \leq b^{-1} \sum_k \int_R |(T_{k/b}f)G_k| =$$

$$b^{-1} \sum_k \int_R |f||\overline{G_{-k}}| \leq b^{-1} \int_R |f| \sum_k |G_k| \leq b^{-1}B \int_R |f| = b^{-1}B\|f\|_{L^1}.$$

This shows that $S$ is a bounded linear operator from a dense subspace of $L^1$ to itself.

*Case II: The $L_0^\infty$-Case.*

For any compactly supported $f \in L_0^\infty$ we have

$$\|Sf\|_{L^\infty} = \|Lf\|_{L^\infty} = \|\sum_k (T_{k/b}f)G_k\|_{L^\infty} = \text{ess sup}|\sum_k (T_{k/b}f)G_k|$$

$$\leq \text{ess sup} \sum_k |T_{k/b}f||G_k| \leq \text{ess sup} \sum_k \|f\|_{L^\infty}|G_k|$$

$$\leq \|f\|_{L^\infty} \left( \text{ess sup} \sum_k |G_k| \right) \leq B\|f\|_{L^\infty}.$$

This makes $S$ a bounded linear operator from a dense subspace of $L_0^\infty$ to itself.

$(2) \Rightarrow (1)$: We assume that $S$ is a bounded linear operator on $L_0^\infty$. Fix n, let $I = [0, a]$ and choose functions $(f_k)_{k=-n}^n$ satisfying:

(1) $|f_k| = \chi_{[k/b,(k/b)+a]}$,

(2) $(T_{-k/b}f_k)G_k = \chi_{[0,a]}|G_k|$.

Letting $f = \sum_k f_k$, since $ab \leq 1$ we have that $\|f\|_{L^\infty} = 1$. Now we have

$$\|S\| \geq \|Sf\|_{L^\infty} \geq \|\chi_{[0,a]}Sf\|_{L^\infty} = \|\sum_{k=-n}^n (T_{-k/b}f_k)G_k\|_{L^\infty}$$

$$= \|\sum_{k=-n}^n |G_k|\|_{L^\infty} = \text{ess sup} \sum_{k=-n}^n |G_k|.$$

This is all we need to establish (1).

For the moreover part of the theorem we apply the theorem to $S^{-1}$ which also satisfies the uniform CC-condition by [39], Theorem 4.14.     $\square$

There is also a classification for the WH-frames $(g, a, b)$ for which the frame operator extends to a bounded linear operator on the Wiener amalgam space [39], Theorem 7.2.

**Theorem 9.19.** *If $ab \leq 1$ and $g \in \mathbf{PF}$, the following are equivalent:*

*(1) The frame operator $S$ is a bounded linear operator from $W(L^\infty, \ell_1)$ to itself.*

*(2) We have*

$$\sum_{k \in \mathbb{Z}} \|G_k\|_\infty = \sum_{k \in \mathbb{Z}} ess \ sup |G_k(t)| = B < \infty.$$

Now let us return to alternate dual frames as given in Definition 4.17, but for the case of Weyl-Heisenberg frames. We start with an important result of Wexler-Raz [135]. The rigorous form of the result given below is due to Janssen [101].

**Theorem 9.20.** *Let $g, h \in \mathbf{PF}$. Then $(E_{mb}T_{na}h)$ and $(E_{mb}T_{na}g)$ are alternate dual WH-frames if and only if $h \perp E_{n/a}T_{m/b}g$, for all $(m, n) \neq (0, 0)$ and $< h, g >= ab$.*

Now we are ready for the classification of alternate dual WH-frames. Parts of the following result were done in an even more general setting by Ron and Shen [127] and Janssen [74].

**Theorem 9.21.** *For $(g, a, b)$ a WH-frame with frame operator $S$ and $h \in \mathbf{PF}$ the following are equivalent:*

*(1) $(E_{mb}T_{na}h)$ is an alternate dual frame for $(E_{mb}T_{na}g)$.*

*(2) We have*

    *(a) $\sum_{n \in \mathbb{Z}} h(x - na)\overline{g(x - na)} = b$ a.e.*

    *(b) $\sum_{n \in \mathbb{Z}} h(x - na)\overline{g(x - na - k/b)} = 0$ a.e. for all $k \neq 0$.*

*(3) $h = S^{-1}g + f$, where $f \in L^2(\mathbb{R})$ and $f \perp span_{n,m \in \mathbb{Z}} E_{n/a}T_{m/b}g$.*

It can be shown that the canonical dual generated by $S^{-1}g$ has the least $L^2$-norm among all $h \in L^2(\mathbb{R})$ satisfying $h \perp E_{n/a}T_{m/b}g$, $(n, m) \neq (0, 0)$, and $< h, g >= ab$. For proofs of this fact we refer to [102], Proposition 3.3 and [56], Proposition 4.2. Note that $(E_{mb}T_{na}g)$ is a normalized tight frame if and only if we can replace $h$ in Theorem 9.21 by the function $g$. Also in this case $S = I$ so part (3) of the theorem becomes: $h = g + f$ where $f \perp span_{n,m \in \mathbb{Z}} E_{n/a}T_{m/b}g$. As we have seen, the canonical dual to a WH-frame is again a WH-frame, but Li [108] has shown there exist alternate duals not having the WH-structure.

We mention again that an important question in WH-frame theory is to find $S^{-1}g$. Strohmer (see Chapter 8 of [74]) has developed a series of algorithms for carrying out this computation. Again, this is an important and fruitful area of study. Also in [74] one can find a theory of multi-window WH-frames as well as applications to optics, signal detection and a host of other directions.

This ends our tour of Weyl-Heisenberg frame theory. Admittedly, we have left out the most important part of the theory: applications to image processing, signal processing, data compression etc. Unfortunately, these topics would require significantly more space than we have available here, and because of the rather specific nature of each set of applications we leave it to the reader to consult the myriad of books devoted to each one.

## 10. Topics in Modern Abstract Frame Theory

It was an open question for a time whether every bounded frame must contain (i.e. have a subsequence which is) a Riesz basis. Seip [130] proved there are exponential frames for $L^2[\pi, \pi]$ which do not contain a Riesz basis. Casazza and Christensen [33, 35] gave elementary examples of this phenomenon.

**Example 10.1.** *The sequence of vectors* $((f_j^n)_{j=1}^{n+1})_{n=1}^\infty$ *(see equations (4.1) for the definition) forms a normalized tight frame for* $H = (\sum_n \oplus H_n)_{\ell_2}$ *which does not contain a Schauder basis.*

A more delicate example [28] is that of a normalized tight frame which contains a Schauder basis but does not contain a Riesz basis. This means that our frame has subsequences which are $\omega$-independent, but any subset of the frame which is $\omega$-independent, is no longer a frame for $H$.

One class of frames which always contain a Riesz basis are the **Besselian frames** introduced by Holub [98]. A frame $(f_i)$ is **Besselian** if whenever $\sum_i a_i f_i$ converges, then $(a_i) \in \ell_2$. We have [98]

**Proposition 10.2.** *If* $(f_i)$ *is a frame for* $H$, *the following are equivalent:*

*(1)* $(f_i)$ *is a near-Riesz basis (i.e. After deleting a finite number of elements, $(f_i)$ becomes a Riesz basis for $H$).*

*(2)* $(f_i)$ *is Besselian.*

*(3)* $\sum_i a_i f_i$ *converges in $H$ if and only if $(a_i) \in \ell_2$.*

There is room here for some interesting work still to be done.

**Problem 10.3.** *Classify the Weyl-Heisenberg frames which contain a Riesz basis.*

We should point out that Problem 10.3 may be exceptionally difficult beyond its usefulness to the theory. However, it is possible that some new important properties of WH-frames will show up in its solution.

**Problem 10.4.** *Classify the exponential frames (or frames of translates) which contain Riesz bases for their span.*

In [46], Christensen introduces the notion of a Riesz frame. We call $(f_i)$ a **Riesz frame** if $(f_i)$ is a frame with frame bounds $A, B$ and every subset of $(f_i)$ is a frame for its closed linear span with the same frame bounds $A, B$. A slight weakening of this definition is: We say that a frame $(f_i)$ has the **subframe property** if every subset of $(f_i)$ is a frame for its closed linear span. Christensen [46] showed that every Riesz frame contains a Riesz basis.

A simple proof is just to choose a maximal linearly independent subset of a Riesz frame, and note that it follows from Proposition 4.3 that this subset is a Riesz basis. It was shown in [30] that every frame with the subframe property also contains a Riesz basis. In fact, a complete classification of frames with the subframe property (but using Riesz frames to do it) is in [27]. However, we do not have a good classification of Riesz frames.

**Problem 10.5.** *Classify the Riesz frames. Are there exponential frames with this property (which do not degenerate into Riesz bases)?*

As we have seen, one of the main difficulties with frame theory is that in order to do *reconstruction* we have to invert the frame operator $S$. This is often difficult if not impossible. It could be quite useful if this inversion could be done by finite dimensional methods. This is the idea behind the **projection methods** introduced by Christensen [44].

We say that a frame $(f_i)$, with frame operator $S$, satisfies the **projection method** if

$$\lim_{n \to \infty} < f, S_n^{-1} f_i > = < f, S^{-1} f >, \text{ for all } i, \ f \in H.$$

where $S_n$ is the frame operator on $K_n = \text{span}(f_i)_{i=1}^n$. We say that $(f_i)$ satisfies the **strong projection method** if there is an increasing sequence of subsets of $\mathbb{N}$, $I_1 \subset I_2 \subset \cdots \nearrow \mathbb{N}$, and frame operators $S_n$ for $K_n = \text{span}_{i \in I_n} f_i$ (i.e. $S_n f = \sum_{i \in I_n} < f, f_i > f_i >$) so that

$$\sum_{i \in I_n} | < f, S_n^{-1} f_i - < f, S^{-1} f_i > |^2 \to 0, \text{ for all } f \in H.$$

We call $(f_i)$ a **conditional Riesz frame** if the frames $(f_i)_{i \in I_n}$ have common frame bounds. This is equivalent to saying:

$$\sup_n \|S_n^{-1}\| < \infty.$$

Christensen [44] shows that $(f_i)$ satisfies the projection method if and only if there are constants $C_j$ so that

$$\|S_n^{-1} f_j\| \le C_j, \text{ for all } n \text{ such that } j \in I_n.$$

These ideas are related to the solution of **moment problems**. That is, given $(a_i) \in \ell_2$, when does there exist an $f \in H$ so that

$$< f, f_i > = a_i, \text{ for all } i \in I?$$

Some deeper results on moment problems can be found in [129]. We combine several results from [34] below.

**Theorem 10.6.** *Let $(f_i)$ be a frame for a Hilbert space $H$. The following are equivalent:*

*(1) $(f_i)$ is a conditional Riesz frame.*

*(2) The strong projection method works.*

*(3) $S^{-1} P_n f \to S^{-1} f$, for all $f \in H$ (where $P_n$ is the orthogonal projection of $H$ onto $span_{i \in I_n} f_i$).*

*(4) $< S_n^{-1} P_n f, g > \to < S^{-1} f, g >$, for all $f, g \in H$.*

*(5) $\sum_{i \in I_n} a_i S_n^{-1} f_i \to \sum_{i \in I} a_i S^{-1} f_i$, as $n \to \infty$ and for all $(a_i) \in \ell_2$.*

*(6) $S_n^{-1} \sum_{i \in I_n} b_i f_i \to 0$, as $n \to \infty$ for all $(a_i) \in \ell_2$ for which $\sum_i a_i f_i = 0$.*

There are several problems with the projection methods. For one, there appear to be only a few specialized frames which satisfy them. It is not known, except in trivial cases, if Weyl-Heisenberg frames have this property. Variations of this method which work for any frame appeared in [37]. However, this is an area which needs much more attention. The importance here is that these results are in the direction of finite frame theory (which is really all there is in the concrete applications).

There are several more abstract directions of research for frame theory. One of them is in the direction of what is called "local theory" in Banach space theory. That is, frames for finite dimensional spaces. As we have seen, there are normalized frames for a Hilbert space which do not contain a Riesz basis. In the finite dimensional setting, such questions become "quantitative" in nature. That is, any frame for a finite dimensional Hilbert space *always* contains a Riesz basis - namely any subset which is complete and linearly independent. So here we are really asking: If $(f_i)$ is a frame for a finite dimensional Hilbert space $H$, is there a subset of $(f_i)$ which is a Riesz basis for $H$ for which the Riesz basis constants are a function of the frame bounds of $(f_i)$ but independent of $(f_i)$ and the dimension of the Hilbert space? It is easily seen that this problem also has a negative answer [29]. But, using some deep results from the local theory of Banach spaces, Casazza [29] showed that for normalized frames $(f_i)$ for $n$-dimensional Hilbert space $H_n$ and all $\epsilon > 0$, there is a subset of $(f_i)$ which is a Riesz basis for a subspace $K_n$ of $H_n$ with $\dim K_n \geq (1 - \epsilon)n$ and having Riesz basis constants only a function of the frame bounds. Vershynin [133] removed the unnecessary hypothesis that the frame be normalized in the above theorem. This produced a powerful method for working with the "John's ellipsoid" in Banach space theory and produced some significant simplifications of several formally quite deep results in this area. In general, however, frame theory for finite dimensional Hilbert spaces is lacking in deep quantitative estimates which could be quite useful - especially in WH-frame theory.

Frames for Hilbert spaces were generalized to Banach spaces by Gröchenig [87]. A **Banach space of scalar valued sequences** (often called a **BK-space**) is a linear space of sequences with a norm which makes it a Banach space and for which the coordinate functionals are continuous.

**Definition 10.7.** *Let $X$ be a Banach space and let $X_d$ be an associated Banach space of scalar valued sequences indexed by $\mathbb{N}$. Let $(y_i)$ be a sequence of elements from $X^*$ and $(x_i)$ be elements of $X$. If*

*(1) $\{(x, y_i)\} \in X_d$, for each $x \in X$,*

*(2) The norms $\|x\|_X$ and $\|\{(x, y_i)\}\|_{X_d}$ are equivalent. i.e. there are constants $A, B > 0$ so that for all $x \in X$ we have*

$$A\|x\|_X \leq \|\{(x, y_i)\}\|_{X_d} \leq B\|x\|_X.$$

*(3) For all $x \in X$ we have $x = \sum_i < x, y_i > x_i$,*

*then we call $((y_i), (x_i))$ an* **atomic decomposition** *of $X$ with respect to $X_d$.*

Gröchenig [87] also defines the notion of a **Banach frame**. Before smooth wavelet orthonormal bases were discovered for Hilbert space $L^2(\mathbb{R})$, Frazier and Jawerth [78] constructed wavelet atomic decompositions for Besov spaces, which they called the $\phi$-transform. A host of papers have appeared surrounding these notions [47, 49, 66, 134]. Also, Feichtinger and Gröchenig developed a quite general theory for a large class of functions and group representations [69, 70, 71]. (See also Chapters 3 and 5 of [74]). This has proved to be quite useful in the study of Besov spaces. This was important since Walnut [134] showed that there are no Gabor frames for weighted Hilbert space. Casazza, Han and Larson [41] made a detailed study of a host of possible generalizations to Banach space theory of the notion of a frame including atomic decompositions.

Finally, we mention that a group of strong people in operator theory have started to enter frame theory bringing with them some very powerful tools from operator theory - especially $C^*$-algebra and Von Neuman algebra theory (See for example, Frank and Larson [76], Frank, Paulsen and Tiballi [77], Han and Larson [93], and Gabardo and Han [82], and results of Z. Landau, R. Balan etc). It is not clear where this will go yet, but it is clear that this development will produce some very deep new results in the area.

DEPARTMENT OF MATHEMATICS, UNIVERSITY OF MISSOURI-COLUMBIA, COLUMBIA, MO 65211

*E-mail address*: pete@math.missouri.edu